\documentclass[a4paper,11pt,reqno]{amsart}
\usepackage[centering,includeheadfoot,margin=2.5 cm]{geometry}
\usepackage[percent]{overpic}
\usepackage{todonotes}
\usepackage{enumitem}
\usepackage{graphics,enumitem,epsfig,textcomp}
\usepackage{amsfonts,amsmath,amssymb,euscript,color,mathrsfs, bbm}
\usepackage[utf8]{inputenc}	
\usepackage{float} 

\usepackage[caption=false,subrefformat=parens,labelformat=parens]{subfig}
\usepackage{url}
\usepackage[symbol]{footmisc}
\usepackage{hyperref}
\usepackage{cases}
\usepackage{bigints}
\hypersetup{
    colorlinks=true,
    linkcolor=blue,
    filecolor=blue,      
    urlcolor=blue,
    citecolor=blue,
    }
\usepackage{comment}

\usepackage[square,numbers]{natbib}
\usepackage{amsmath,amsthm,amssymb,amsfonts,lineno}
\usepackage{esint}
\usepackage{braket, leftindex}

\numberwithin{equation}{section}
\newtheorem{theorem}{Theorem}

\newtheorem{lemma}{Lemma}
\newtheorem{proposition}{Proposition}
\newtheorem{remark}{Remark}

\def\diver{\mathrm{div}}

\def\d{\,\mathrm{d}}
\def\eps{\varepsilon}
\def\N{\mathbb{N}}
\def\R{\mathbb{R}}
\def\C{\hbox{\rlap{\kern.24em\raise.1ex\hbox{\vrule height1.3ex width.9pt}}C}}

\def\P{\hbox{\rlap{I}\kern.16em P}}
\def\Q{\hbox{\rlap{\kern.24em\raise.1ex\hbox
      {\vrule height1.3ex width.9pt}}Q}}
\def\M{\hbox{\rlap{I}\kern.16em\rlap{I}M}}
\def\Z{\hbox{\rlap{Z}\kern.20em Z}}
\def\({\begin{eqnarray}}
\def\){\end{eqnarray}}
\def\[{\begin{eqnarray*}}
\def\]{\end{eqnarray*}}
\def\part#1#2{\frac{\partial #1}{\partial #2}}

\def\grad{\nabla}

\def\pmb#1{\setbox0=\hbox{$#1$}
  \kern-.025em\copy0\kern-\wd0
  \kern-.05em\copy0\kern-\wd0
  \kern-.025em\raise.0433em\box0 }
\def\bar{\overline}

\def\weakconv{\rightharpoonup}

\def\d{\,\mathrm{d}}
\def\N{\mathbb{N}}
\def\R{\mathbb{R}}

\def\epsilon{\varepsilon}

\def\E{\mathcal{E}}
\def\P{\mathbb{P}}
\def\Q{\mathbb{Q}}

\usepackage{cancel}
\title{Measure-based approach to mesoscopic modeling of optimal transportation networks}
\date{January 2024}

\begin{document}
\pagenumbering{gobble}
\maketitle
\pagenumbering{arabic}

\centerline{
     {\large Jan Haskovec}\footnote{Mathematical and Computer Sciences
            and Engineering Division,
         King Abdullah University of Science and Technology,
         Thuwal 23955-6900, Kingdom of Saudi Arabia;
         {\it jan.haskovec@kaust.edu.sa}}\qquad
     {\large Peter Markowich}\footnote{Mathematical and Computer Sciences and Engineering Division,
         King Abdullah University of Science and Technology,
         Thuwal 23955-6900, Kingdom of Saudi Arabia;
         {\it peter.markowich@kaust.edu.sa}, and
         Faculty of Mathematics, University of Vienna,
        Oskar-Morgenstern-Platz 1, 1090 Vienna;
         {\it peter.markowich@univie.ac.at}}\qquad
         {\large Simone Portaro}\footnote{Mathematical and Computer Sciences
            and Engineering Division,
         King Abdullah University of Science and Technology,
         Thuwal 23955-6900, Kingdom of Saudi Arabia;
         {\it simone.portaro@kaust.edu.sa}}
     }
\medskip \medskip

\textbf{Abstract.}
We propose a mesoscopic modeling framework for optimal transportation networks with biological applications.
The network is described in terms of a joint probability measure on the phase space of tensor-valued conductivity and position in physical space.
The energy expenditure of the network is given by a functional consisting of a pumping (kinetic) and metabolic power-law term, constrained by a Poisson equation accounting for local mass conservation.
We establish convexity and lower semicontinuity of the functional on appropriate sets.
We then derive its gradient flow with respect to the 2-Wasserstein topology on the space of probability measures, which leads to a transport equation, coupled to the Poisson equation.
To lessen the mathematical complexity of the problem, we derive a reduced Wasserstein gradient flow, taken with respect to the tensor-valued conductivity variable only. We then construct equilibrium measures of the resulting PDE system.
Finally, we derive the gradient flow of the constrained energy functional with respect to the Fisher-Rao (or Hellinger-Kakutani) metric, which gives a reaction-type PDE. We calculate its equilibrium states, represented by measures concentrated on a hypersurface in the phase space.
\vspace{3mm}

\section{Introduction}\label{sec:intro}
In this paper we propose a mesoscopic model describing optimal transportation networks, primarily in biological context - leaf venation in plants, blood circulatory systems in animals or neuronal networks.
The model is based on the discrete (microscopic) framework introduced by Hu and Cai in \cite{hucai2013}. They considered an energy functional posed on a discrete graph, consisting of a pumping and a metabolic term. The functional is coupled to the Kirchhoff law imposing global mass conservation.
The model has been further studied in \cite{kreusser2019, Hu-Cai-19}. In \cite{burger2019} it has been shown that the local energy minimizers in the case of concave metabolic cost are trees (i.e., loop-free graphs), and vice versa, that every tree induces a local minimizer \cite{HV23}.

A continuum version of the model has been introduced and thoroughly studied, analytically and numerically, in the series of papers \cite{albi2000, astuto1, astuto-port, Kreusser-rigorous, HMP15, HMPS16}. Here the network is understood as a porous medium with principal permeability given by a vector field on a spatial domain. The continuum version of the energy functional is constrained by a Poisson equation for the material pressure.
A more general model has been derived
in \cite{CMS-2022, portaro2023}, where the permeability of the porous medium is given by a symmetric positive semidefinite tensor field.

In \cite{burger2019} a mesoscopic description has been proposed, modeling the network in terms of a probability that at a given location the network has a vector-valued direction and a scalar nonnegative conductivity. Consequently, the system is described by a time-dependent probability measure posed on the respective phase space and its evolution is subject to a hyperbolic transport equation. This is again coupled to a Poisson equation with measure-valued permeability tensor, obtained as a moment of the probability measure.
In \cite{burger2019}, formal connection to the discrete model has been established and its stationary solutions provided, including discrete network solutions.

In this paper we propose a fully general mesoscopic model where the permeability of the medium is described by a symmetric positive semidefinite tensor $\C\in\R^{d\times d}$, where $d\geq 1$ is the physical space dimension. Consequently, the network is modeled in terms of a probability measure $\mu=\mu(x,\C)$, where $x\in\bar\Omega$ with $\Omega\subset\R^d$ is a bounded domain.
This represents the joint probability of the network to have conductivity $\C$ at the point $x \in \overline\Omega$. The marginal probability $\mu_{\C}(x) := \int \mu (x, \d \C)$ describes then the position density of the network in $\overline{\Omega}$ whose variations may be caused by different material properties. 
This probabilistic interpretation is analogous to the conventional kinetic interpretation of gas dynamics in which (particle) mass density is calculated by integrating the joint probability with respect to the velocity variable, see, e.g., \cite{Cercignani}. Thus, drawing an analogy, we interpret $\mu$ as a joint probability measure describing local network densities and conductivities.

We first formulate the mesoscopic analogue of the discrete energy functional, constrained by a Poisson equation accounting for local mass conservation. The permeability tensor in the Poisson equation is the first-order $\C$-moment of the measure $\mu$. We establish essential mathematical properties of the energy functional, in particular, its convexity and lower semicontinuity on appropriate sets.
We then calculate the gradient flow of the constrained problem with respect to the 2-Wasserstein topology on the space of probability measures.
The resulting system of partial differential equations can be seen as a model that accounts for dynamic changes in the transportation network.

The choice of the 2-Wasserstein gradient flow is based on the fundamental modeling paradigm to emphasize the aspects of (mass) transportation and associated diffusion processes in the network. 
Moreover, compared to the $p$-Wasserstein distance with a general exponent $p\in[1,\infty]$, the choice $p=2$ facilitates
relative ease of calculations, benefiting from an established analytical framework \cite{ambrosiogigli, ambrosio_handbook, villani2021}.

One of the main difficulties for analysis of gradient flow with the 2-Wasserstein gradient flow of the energy is the presence of
third derivatives of the material pressure in the equation for $\mu$, with nonlinear second order terms under the divergence. This is the main obstacle for developing a well-posedness theory for the corresponding PDE system.
Moreover, as we shall show, no generic steady states exist.
Therefore, in order to lessen the mathematical complexity of the problem, we consider a reduced Wasserstein gradient flow, taken with respect to the $\C$-variable only, while the position $x$ is treated as a parameter.
Then, for every admissible network position density we are able to find a unique equilibrium measure (of monokinetic form, concentrated on a $d$-dimensional manifold in $d+d(d+1)/2$-dimensional phase space) of the resulting PDE system, where the equilibrium pressure is obtained as solution of a $p$-Laplace-type equation.
Moreover, we show that the model of \cite{burger2019} is a special case of the reduced model.

Finally, we derive the gradient flow
of the constrained energy functional with respect to the Fisher-Rao (or Hellinger-Kakutani) metric.
This is the unique (up to a multiplicative factor) smooth metric
which is invariant under the action of the diffeomorphism group on manifolds of dimension higher than one, see \cite{bauer}.
Compared to the 2-Wasserstein metric, which induces a transport-type PDE as the gradient flow,
the Fisher-Rao prioritizes reactive processes to drive the network transportation dynamics. Indeed, the corresponding gradient flow equation is of reaction-type. Also here we benefit from an existing toolbox for calculus and analysis, see, e.g., \cite{liero2016, liero2018, mielke}. Moreover, let us note that linear combinations of the 2-Wasserstein and Fisher-Rao metric can be used to derive dynamical models, although we do not follow this approach here.

We focus on the equilibrium states of the corresponding gradient flow, represented by measures concentrated on a hypersurface in the phase space. 
In particular, the steady state manifold is parametrized by 
a probability measure on the hypersurface
and a real constant.
This can be seen as one of the most significant results of our paper, beside laying ground for future analytical work (geodesic convexity, well-posedness of the gradient flow equations etc).

The complexity of the set of equilibrium measures shows the richness of our transportation energy model and its power in representing a multitude of network patterns.
Let us note that the measure-based energy functional and the network dynamics resulting from the gradient flows are very costly for numerical treatment due to the high dimensionality of the underlying phase space.
However, simpler and more practically useful models can be easily obtained by restricting the energy functional to appropriate subsets or submanifolds of the probability space. In particular,
the mesoscopic description of \cite{burger2019}, where the network conductivity is vector-valued, which reduces the dimension of the phase space. Also the tensor-based macroscopic description of \cite{CMS-2022} and vector-based models \cite{albi2000, HMP15, HMPS16} can be derived from our model. Consequently, the measure-based approach presented in this paper can be seen as an unifying and generalizing framework for a whole hierarchy of models for optimal transportation networks.

This paper is organized as follows.
In Section \ref{sec:discrete} we describe the discrete model of \cite{hucai2013}.
In Section \ref{sec:meso} we introduce the mesoscopic model, study its mathematical propeties and derive its full and reduced Wasserstein gradient flow.
In Section \ref{sec:stationary} we investigate the stationary solutions of the respective evolution equations.
Finally, in Section \ref{sec:fisher_rao} we derive the Fisher-Rao gradient flow of the energy functional and study its stationary states.

\section{Discrete model}
\label{sec:discrete}
In this section we describe the discrete (microscopic) model introduced by Hu and Cai in \cite{hucai2013} and further studied in \cite{albi2000, burger2019, kreusser2019}.

Let $\mathcal{G} = (\mathcal{V}, \mathcal{E})$ denote an undirected connected graph \cite{gross2013} with a finite set $\mathcal{V}$ of vertices $i \in \mathcal{V}$, and a finite set $\mathcal{E}\subseteq\mathcal{V}\times\mathcal{V}$ of edges, where each edge $(i,j) \in \mathcal{E}$ represents a connection between vertices $i$ and $j$. We assume that each pair of vertices is connected by at most one edge, and we exclude self-loops, i.e., edges of the form $(i,i)$. As the graph is undirected, we refer by $(i,j)$ and $(j,i)$ to the same edge.

To each edge $(i, j)\in\mathcal{E}$ we assign a length parameter $L_{ij} > 0$ (in preparation of the embedding of the graph into a d-dimensional spatial domain), and by $C_{ij} \geq 0$ we denote its conductivity. At each vertex $i \in \mathcal{V}$ we prescribe the strength of the external material flow $S_i\in\R$, where we adopt the convention that $S_i$ is positive for inflow nodes (sources) and negative for sinks. Moreover, for each $i\in\mathcal{V}$ we introduce the material pressure of the flow $P_i \in \R$.

The oriented flux from vertex $i$ to vertex $j$ is denoted by $Q_{ij}$, with $Q_{ij} = -Q_{ji}$. We assume that the flow takes place in the laminar (Poiseuille) regime characterized by low Reynolds numbers.
Then, the flux through an edge is proportional to the edge conductivity and the pressure drop between its endpoints,
\begin{align}
    \label{eq:discrete_flow}
    Q_{ij} := C_{ij} \frac{P_j - P_i}{L_{ij}} \quad \forall (i, j) \in \mathcal{E}.
\end{align}
At each vertex $i \in \mathcal{V}$ we impose the local conservation of mass, expressed in terms of Kirchhoff Law
\begin{align}
    \label{eq:kirchhoff_Law}
    - \frac{1}{L_{ij}} \sum_{j \in \mathcal{N}(i)} Q_{ij} = - \frac{1}{L_{ij}} \sum_{j \in \mathcal{N}(i)} C_{ij} \frac{P_j - P_i}{L_{ij}} = S_i \qquad \forall i \in \mathcal{V},
\end{align}
where $\mathcal{N}(i)$ denotes the set of vertices connected to $i$ by an edge, $\mathcal{N}(i):=\{j\in\mathcal{V}:(i,j)\in\mathcal{E}\}$. Clearly, a necessary (but not sufficient) condition for the solvability of \eqref{eq:kirchhoff_Law} is the global conservation of mass
\begin{align}
    \label{eq:global_conservation}
    \sum_{i \in \mathcal{V}} S_i = 0,
\end{align}
which we assume in the sequel. For a given vector of conductivities $C = (C_{ij})_{(i,j) \in \mathcal{E}}$, the linear system \eqref{eq:kirchhoff_Law} possesses a solution $(p_i)_{i \in \mathcal{V}}$ if the underlying graph is connected, taking into account only edges with positive conductivities $C_{ij}$, i.e., disregarding edges with zero conductivities. If it exists, the solution $(p_i)_{i \in \mathcal{V}}$ is unique up to an additive constant.

Hu and Cai \cite{hucai2013} proposed the following functional for the energy expenditure of the transportation network,
\begin{align}
    \label{eq:discrete_energy}
    E[C] = \sum_{(i,j) \in \mathcal{E}, i<j} \left(  \frac{Q_{ij}^2[C]}{C_{ij}} + \frac{\nu}{\gamma} C_{ij}^{\gamma} \right) L_{ij},
\end{align}
where $i < j$ in the summation symbol means that each edge in the graph is counted exactly once.
The first term in the expression corresponds to the energy needed to pump the material through an edge $(i,j) \in \mathcal{E}$, i.e., the kinetic energy of the flow through the edge. By the Joule's law, it is given by the product of the pressure drop $(P_j - P_i)$ across the edge and the flow rate $Q_{ij}$,
\begin{align*}
    (P_j - P_i) Q_{ij} = \frac{Q_{ij}^2}{C_{ij}} L_{ij}.
\end{align*}
The second term in the energy incorporates the metabolic cost to sustain the edge, with $\nu > 0$ denoting the metabolic coefficient. The metabolic cost is assumed to be proportional to the edge length $L_{ij}$ and the conductivity $C_{ij}$ raised to the power $\gamma>0$.
The value of $\gamma$ depends on the specific biological application. For instance, in blood vessel systems, $\gamma = \frac{1}{2}$ due to the metabolic cost being proportional to the vessel's cross-section area, see \cite{murray1926} for further details.

Let us note that the gradient flow of the energy \eqref{eq:discrete_energy}, constrained by the Kirchhoff's law \eqref{eq:kirchhoff_Law}, can be explicitly calculated \cite{kreusser2019}. It is given by the following ODE system for the conductivities $C_{ij}=C_{ij}(t)$,
\begin{align}
\label{eq:discrete_GF}
    \frac{d C_{ij}}{dt} = \left( \frac{Q_{ij}[C_{ij}]^2}{ \left( C_{ij} \right)^2} - \nu C_{ij}^{\gamma-1} \right) L_{ij}.
\end{align}

\section{Mesoscopic model and Wasserstein gradient flow}\label{sec:meso}

Now, we introduce the mesoscopic analogue of the discrete model discussed in Section \ref{sec:discrete}. Our approach is inspired by the kinetic model studied in \cite{burger2019}. However, in this work we embed the graph $\mathcal{G}$ into $\overline{\Omega} \subset \R^d$ where $d$ is the spatial dimension and describe the local network conductivity in terms of a symmetric, positive semidefinite tensor $\C\in\R^{d\times d}$. This is significantly more general than \cite{burger2019}, where the conductivity is restricted to the set of rank-one tensors.
We shall show that the model of \cite{burger2019} can be obtained as a particular case of our approach.

The main idea lies in re-interpreting the discrete pumping term through the edge $(i,j)$ with endpoints located in $x_i$, $x_j\in\R^d$ as
\begin{align*}
    \frac{Q_{ij}^2}{C_{ij}} = C_{ij} \frac{\left( P_j - P_i \right)^2}{L_{ij}^2} \approx \grad p \cdot C_{ij} \frac{(x_i-x_j)\otimes(x_i-x_j)}{|x_i-x_j|^2} \grad p,
\end{align*}
where $p=p(x)$ is the continuum material pressure field in the medium. Note that we obviously have $|x_i-x_j|=L_{ij}$.
This construction was used in \cite{CMS-2022} to formally derive the tensor-valued continuum limit of the model, applying tools from finite element approximation.
Here we make a further step in the generalization of the model and formulate the pumping term as $\grad p\cdot\C\grad p$ with $\C \in \overline{\mathscr{S}^d}(\R)$ the network conductivity tensor,
where $\overline{\mathscr{S}^d}(\R)$
denotes the closed convex set of symmetric positive semidefinite matrices in $\R^{d\times d}$.
We reiterate that the network is contained within a fixed bounded domain $\Omega \subset \R^d$.

We now introduce the probability measure $\mu \in \mathscr{P} := \mathscr{P} \left( \overline{\Omega} \times \overline{\mathscr{S}^d}(\R) \right)$, where $\mu(x, \C)$ represents the joint probability of the network to have conductivity $\C\in\overline{\mathscr{S}^d}(\R)$ at the point $x \in \overline\Omega$ and $\mathscr{P}$ is the space of Borel probability measures on $\overline{\Omega} \times \overline{\mathscr{S}^d}(\R)$. 
The marginal probability, defined formally as $\mu_{\C}(x) := \int_{\overline{\mathscr{S}^d}(\R)} \mu (x, \d \C)$ and rigorously by $\mu_{\C}(x) = (\pi_x)_{\#} \mu$, where $\pi_x : \overline{\Omega} \times \overline{\mathscr{S}^d}(\R) \rightarrow \overline{\Omega}$ is the projection $\pi_x (x, \C) = x$, describes then the position density of the network in $\overline{\Omega}$.
We formulate the mesoscopic analogue of the discrete energy functional \eqref{eq:discrete_energy} as\begin{align}
    \label{eq:energy}
    E[\mu] = \int_{\overline{\mathscr{S}^d}(\R)} \int_{\overline{\Omega}} \left(  \grad p[\mu] \cdot \C \grad p[\mu] + \frac{\nu}{\gamma} | \C |^{\gamma} \right) \d \mu + \int_{\Omega}  r | \grad p [\mu] |^2 \d x,
\end{align}
with $\nu > 0$ the metabolic constant. Here we introduce the constant $r \ge 0$ to represent the isotropic background permeability of the network's medium, see, e.g., \cite{HMP15, HMPS16, CMS-2022, portaro2023}.
The mesoscopic energy functional is constrained by the appropriate continuous reformulation of the Kirchhoff's law \eqref{eq:kirchhoff_Law}, which is given by the Poisson equation
\begin{align}
    \label{eq:poisson1}
    - \diver \left( \left( \P[\mu] + r \mathbb{I} \right) \grad p\right) = S,
\end{align}
where $S \in L^2(\Omega)$ is the distribution of sources and sinks and the permeability tensor $\P[\mu]$ is defined as
\begin{align}
    \label{eq:permeability_tensor}
    \P [\mu](x) = \int_{\overline{\mathscr{S}^d} (\R)} \C\, \mu(x, \d \C),
\end{align}
representing the transport directions of the network. Notice that $\P[\mu]$ is the first moment of $\mu$ with respect to $\C$. To close the Poisson equation, we impose the zero-flux boundary condition
\begin{align}
    \label{eq:bc_poisson}
    n \cdot \left( \P[\mu] + r \mathbb{I} \right) \grad p = 0 \quad \mathrm{on} \, \, \partial \Omega,
\end{align}
where $n$ is the outer unit normal vector to the boundary $\partial \Omega$. Obviously, a necessary condition for solvability of \eqref{eq:poisson1}--\eqref{eq:bc_poisson} is the global
mass conservation
\begin{align}
\label{eq:cont_globalconservation}
    \int_{\Omega} S \d x = 0.
\end{align}
Note that the isotropic background permeability term $r\mathbb{I}$ guarantees uniform ellipticity of \eqref{eq:poisson1}.
Sufficient conditions for solvability of the boundary value problem \eqref{eq:poisson1}--\eqref{eq:bc_poisson}, in particular with regard to the regularity properties of the permeability tensor \eqref{eq:permeability_tensor}, shall be discussed below.

We have $\overline{\mathscr{S}^d}(\R) = \mathscr{S}^d(\R) \cup \partial \mathscr{S}^d(\R)$, where $\mathscr{S}^d(\R)$ represents the set of all symmetric, positive definite matrices, and $\partial \mathscr{S}^d(\R)$ denotes its boundary. The latter consists of singular symmetric, positive semidefinite matrices, i.e., possessing at least one zero eigenvalue. We note that the set $\mathscr{S}^d(\R)$ forms the interior of a convex cone, which is a differentiable manifold of dimension $\frac{d(d+1)}{2}$. This space, equipped with the Euclidean distance constitutes a complete metric space. Furthermore, $\overline{\mathscr{S}^d}(\R)$ is a semialgebraic set, see, e.g., \cite{malick2012}. Semialgebraic sets can be decomposed into a union of connected smooth manifolds, referred to as ``strata". Let ${\mathscr{R}_r}$ denote the smooth submanifold composed of positive semidefinite matrices of rank $r$, where $r = 0, 1, \dots, d$. We then express $\overline{\mathscr{S}^d}(\R)$ as
\begin{align*}
    \overline{\mathscr{S}^d} (\R) = \bigcup_{r=0}^d {\mathscr{R}_r} = \bigcup_{r=0}^{d-1} {\mathscr{R}_r} \cup \mathscr{S}^d (\R).
\end{align*}
This implies that $\partial \mathscr{S}^d(\R) = \bigcup_{r=0}^{d-1} {\mathscr{R}_r}$, signifying that the boundary is a finite union of smooth manifolds. Henceforth, this space enables the application of the divergence theorem, which we will frequently make use of in the sequel.

We shall also use the vectorization of a symmetric matrix. Given a matrix $A \in \R^{d \times d}$, we uniquely map it to a vector, denoted as $\mathsf{vec} (A) \in \R^{d^2}$. 
This map is realized by arranging the elements of $A$ row-wise in $\mathsf{vec} (A)$. However, because we work exclusively with symmetric matrices, we can reduce the dimension of the resulting vector to $\frac{d(d+1)}{2}$ by employing the transformation $\mathsf{v}(A) = D_d^T \mathsf{vec}(A)$ for $A$ symmetric. 
Here, $D_d^T \in \R^{\frac{d(d+1)}{2} \times d^2}$ represents the duplication matrix that removes the redundant elements of $A$ \cite{magnus1980}.

\subsection{Gradient flow in the $2$-Wasserstein topology}\label{subsec:GFWass}

We define the space $X := \overline{\Omega} \times \overline{\mathscr{S}^d}(\R)$ and call $\xi = \left( x, \mathsf{v}(\C) \right) \in X$ a generic element of $X$. We equip $X$ with the Euclidean distance $d(\xi_1, \xi_2)$ for all $\xi_1, \xi_2 \in X$.
We call $\gamma \in \mathscr{P}(X \times X)$ a transport plan between $\mu, \nu \in \mathscr{P}(X)$ if
\begin{align*}
    \int_{X \times X} \varphi (\xi_1) \d \gamma (\xi_1, \xi_2) = \int_{X} \varphi(\xi_1) \d \mu(\xi_1) \quad \forall \, \varphi : X \rightarrow \R \, \, \mathrm{Borel} \, \, \mathrm{and} \, \, \mathrm{bounded},
\end{align*}
and
\begin{align*}
    \int_{X \times X} \psi (\xi_2) \d \gamma (\xi_1, \xi_2) = \int_{X} \psi(\xi_2) \d \nu(\xi_2) \quad \forall \, \psi : X \rightarrow \R \, \, \mathrm{Borel} \, \, \mathrm{and} \, \, \mathrm{bounded}.
\end{align*}
In other words, $\mu$ and $\nu$ are the marginal probabilities of the transport plan $\gamma$.
We denote by $\Gamma(\mu, \nu)$ the set of transport plans of $\mu$ and $\nu$. Note that $\Gamma(\mu, \nu)$ is always nonempty, since the product measure $\gamma = \mu \otimes \nu$, defined via
\begin{align*}
    \int_{X \times X} \phi(\xi_1, \xi_2) \d \gamma (\xi_1, \xi_2) = \int_{X \times X} \phi(\xi_1, \xi_2) \d \mu (\xi_1) \d \nu (\xi_2) \quad \forall \, \phi : X \times X \rightarrow \R,
\end{align*}
is obviously a transport plan.

Let us define the set of probability measures with finite second-order $\C$-moment
\begin{align*}
    \mathscr{P}_2 (X) := \big{\{} \mu \in \mathscr{P}(X) : \int_X | \C |^2 \ \d \mu < \infty \big{\}}.
\end{align*}
Given measures $\mu, \nu \in \mathscr{P}_2(X)$, we define their $2$-Wasserstein distance as
\begin{align}
    \label{eq:wass_distance}
    W_2 (\mu, \nu) := \left( \inf_{\gamma \in \Gamma(\mu, \nu)} \int_{X \times X} d^2 (\xi_1, \xi_2) \d \gamma (\xi_1, \xi_2) \right)^{\frac{1}{2}}.
\end{align}
The distance $W_2$ induces a topology on $\mathscr{P}_2(X)$, which we call Wasserstein topology. Moreover, we refer to $\left( \mathscr{P}_2 (X), W_2 \right)$ as the $2$-Wasserstein space. For a detailed analysis of the properties of Wasserstein spaces, we refer to \cite{ambrosio_handbook}.


Note that $\mu \in \mathscr{P}_2(X) $ implies that $\P[\mu] \in \mathcal{M}^+ (\overline\Omega)$, which is clearly not sufficient to define a proper notion of solution of the Poisson equation \eqref{eq:poisson1}. Therefore, we redefine the energy functional \eqref{eq:energy} as $E[\mu] := + \infty$ whenever $\P[\mu] \not\in L^1 (\Omega)$,
implying that its convex domain is given by
\begin{align}
    \label{eq:domain_E}
    D(E) := \left\{ \mu \in \mathscr{P}_2(X): \, \P[\mu] \in L^1 (\Omega) \;\mathrm{and}\; \int_{X} |\C|^{\gamma} \d\mu < +\infty \right\}.
\end{align}

Note that every $\C^k$ moment of $\mu \in D(E)$ is absolutely continuous with respect to the Lebesgue measure for $0<k \le \max \{ 2, \gamma \}$ and the position density ($k=0$) is absolutely continuous with respect to the Lebesgue measure on $\Omega$ if and only if the restriction of $\mu$ to $\left( \overline{\Omega}\times\{0\} \right)$ is such.
In fact, adding to $\mu$ any measure which concentrates in $\{ (x, \C) \in X : \C = 0 \}$ does not change $E[\mu]$. Thus we might as well replace the integrability condition on $\mathbb{P}[\mu]$ in $D(E)$ by instead requiring the position density $\mu_{\C}(x) = (\pi_x)_{\#} \mu$ to be in $L^1(\Omega)$.

In the following, we need to establish a well-posedness theorem for the Poisson equation \eqref{eq:poisson1}. To this end we introduce the following subspace of the Sobolev space ${H^1}(\Omega)$,
\begin{align}
    \label{eq:space_H}
    H := \left\{ v \in{H^1}(\Omega): \; \int_{\Omega} \grad v \cdot \P[\mu] \grad v \d x < \infty, \; \int_{\Omega} v \d x = 0 \right\}.
\end{align}
We refer to \cite{ukhlov} for the proof of completeness of $H$, provided that $\P[\mu] \in L^1(\Omega)$.

\begin{lemma}\label{poisson_solvability}
    Let $S \in L^2(\Omega)$, $\P[\mu] \in L^1 (\Omega)$. Then the Poisson equation \eqref{eq:poisson1}-\eqref{eq:bc_poisson} admits a unique weak solution $p \in H$ satisfying
    \begin{align}
        \label{eq:estimate_gradp}
        \| \grad p \|_{L^2(\Omega)} \le \frac{C_{\Omega}}{r} \| S \|_{L^2(\Omega)},
    \end{align}
    where $C_{\Omega} > 0$ is the Poincarè constant.
    \begin{proof}
        Let us define the bilinear form $B : H \times H \rightarrow \R$ as
        \begin{align*}
            B(p, \varphi) := \int_{\Omega} \grad \varphi \cdot \left( \P[\mu] + r \mathbb{I} \right) \grad p \d x.
        \end{align*}
        Consider the space $H$ equipped with the norm induced by $B$, i.e., $\| p \|^2_{H} := B(p, p)$. Coercivity of the bilinear form follows directly from the positive semidefinitness of $\P[\mu]$. To establish continuity, we observe that for all $p, \varphi \in H$
        \begin{align*}
            | B(p, \varphi) | &\le \int_{\Omega} \left( \big{|} \sqrt{\P[\mu]} \grad \varphi \big{|} \big{|} \sqrt{\P[\mu]} \grad p \big{|} + r \big{|} \grad \varphi \big{|} \big{|} \grad p \big{|} \right) \d x \\
            & \le \| p \|_{H} \| \varphi \|_{H},
        \end{align*}
        by the Cauchy-Schwartz inequality.
        Consequently, the result follows by an application of the Lax-Milgram theorem. Furthermore, with the Poincar\'e inequality we have
        \begin{align*}
            \| \grad p \|^2_{L^2(\Omega)} \le \frac{1}{r} B(p, p) = \frac{1}{r} \int_{\Omega} S p \d x \le \frac{C_{\Omega}}{r} \| S \|_{L^2(\Omega)} \| \grad p \|_{L^2(\Omega)},
        \end{align*}
        proving \eqref{eq:estimate_gradp}.
    \end{proof}
\end{lemma}

To establish a formal representation of a gradient flow on the manifold $\mathscr{P}_2(X)$ we follow \cite{mielke}. We define the gradient structure triple $\left(M, E, \mathbb{G} \right)$, where
\begin{align*}
    M := \mathscr{P}_2(X) \cap C^{\infty}_+ \left( X \right),
\end{align*}
with $C^{\infty}_+$ as the space of smooth positive functions, $E$ is the energy functional \eqref{eq:energy}, and $\mathbb{G}(\mu) : T_{\mu} M \rightarrow T_{\mu}^* M$ defines the Riemannian structure (a symmetric and positive linear operator) on $\mathscr{P}_2(X)$. Here $T_{\mu} M$ denotes the tangent space to $M$ at $\mu$, given by $T_{\mu} M = \left\{ \nu : X \rightarrow \R \, \, \mathrm{smooth} \; \mathrm{with} \; \, \int_{X} \d \nu = 0 \right\}$, see \cite{mielke, otto_porous}. The cotangent space is $T^*_{\mu} M = \{ \nu : X \rightarrow \R \; \mathrm{smooth} \} / \R$.
Observe that $\mathbb{G}$ induces a positive, symmetric metric tensor $g_{\mu}$ via
\begin{align}
    g_{\mu}(v, \tilde{v}) := \leftindex_{T_{\mu}^* M} {\braket{ \mathbb{G}(\mu) v, \tilde{v}}}_{T_{\mu} M} \quad \forall \, v, \tilde{v} \in T_{\mu} M, \notag
\end{align}
where $\leftindex_{X^*}{\braket{\cdot, \cdot}}_X$ denotes the duality pairing between $X$ and its dual space $X^*$.
Next, we introduce the Onsager operator $\mathbb{K}(\mu) : D(\mathbb{K}(\mu)) \subseteq T^{*}_{\mu} M \rightarrow T_{\mu} M$ which is the inverse of $\mathbb{G}(\mu)$, i.e., $\mathbb{G}(\mu) = \mathbb{K}(\mu)^{-1}$, whose domain $D(\mathbb{K}(\mu))$ is the range of $\mathbb{G}(\mu)$. 
Otto \cite{otto_porous} used formal Riemannian calculus to show that the $2$-Wasserstein distance $W_2$ is generated when $\mathbb{G}(\mu)$ is the solution operator of the weak formulation of the Neumann problem
\begin{align*}
    \int_{X} \mu \grad_{\left( x, \C \right)} \Phi \cdot \grad_{\left( x, \C \right)} \varphi \d x \d \C = \int_{X} s \varphi \d x \d \C \quad \forall \varphi \in H^1(\Omega) \cap \bigg{\{} \varphi : \int_{X} \varphi \d x \d \C = 0 \bigg{\}}
\end{align*}
where $s \in T_{\mu}M$ and $\Phi = \mathbb{G}(\mu) s \in T^{*}_{\mu} M$. Thus, we can identify (on a formal level)
\begin{align*}
    \mathbb{K}(\mu) \Phi = - \diver_{\left(x, \C \right)} \left( \mu \grad_{(x, \C)} \Phi \right)
\end{align*}
and 
\begin{align*}
    D(\mathbb{K}(\mu)) = \{ \Phi : X \rightarrow \R \; \mathrm{smooth} : \mu \grad_{\left(x, \C\right)} \Phi \cdot n_X = 0 \; \mathrm{on} \; \partial X \} / \R,
\end{align*}
where $n_X$ denotes a unit normal vector orthogonal to $\partial X$.

Then, the above defined gradient structure generates the associated gradient-flow equation (the so called $2$-Wasserstein gradient flow)
\begin{align}
    \label{eq:gradient_structure}
    \partial_t \mu = - \mathbb{K}(\mu) \frac{\delta E[\mu]}{\delta \mu} = \diver_{\left(x, \C\right)} \left ( \mu \grad_{\left(x, \C\right)} \frac{\delta E[\mu]}{\delta \mu} \right),
\end{align}
where $\frac{\delta E[\mu]}{\delta \mu} \in D(\mathbb{K}(\mu)) \subseteq T_{\mu}^* M$ is the first variation of $E$ at $\mu$. Thus \eqref{eq:gradient_structure} is subject to a zero outflow boundary condition
\begin{align*}
    \mu \grad_{\left(x, \C\right)} \frac{\delta E[\mu]}{\delta \mu} \cdot n_X = 0 \quad \mathrm{on} \; \partial X.
\end{align*}

We see that $E$ is dissipated along the flow, that is for $\mu : [0, T] \rightarrow M$ the function $t \mapsto E[\mu(t)]$ is non-increasing,
\begin{align}
    \label{eq:flow_dissipation}
    \frac{\d}{\d t} E[\mu(t)] &= \leftindex_{T_{\mu}^* M} {\left\langle \frac{\delta E[\mu]}{\delta \mu}, \partial_t \mu\right\rangle}_{T_{\mu} M} = - \leftindex_{T_{\mu}^* M} {\left\langle \frac{\delta E[\mu]}{\delta \mu}, \mathbb{K}(\mu) \frac{\delta E[\mu]}{\delta \mu} \right\rangle}_{T_{\mu} M} \notag \\
    &= - \leftindex_{T_{\mu}^* M} {\braket{ \mathbb{G}(\mu) \partial_t \mu, \partial_t \mu }}_{T_{\mu} M} \le 0.
\end{align}

\subsection{Mathematical properties of the energy functional}\label{subsec:properties}
Throughout this paper, the first variation of the energy \eqref{eq:energy} plays a crucial role in describing the network's evolutionary dynamics.

\begin{proposition}\label{first_variation}
    The first variation of the energy functional \eqref{eq:energy} constrained by \eqref{eq:poisson1}-\eqref{eq:bc_poisson} is given by
    \begin{align}
        \label{eq:first_variation}
        \frac{\delta E[\mu]}{\delta \mu} = - \grad p[\mu] \cdot \C \grad p[\mu] + \frac{\nu}{\gamma} | \C |^{\gamma}, \; \; \mu \in D(E).
    \end{align}

    \begin{proof}
        Let $\eps \in \R$, $\mu = \mu_0 + \epsilon \mu_1$ and $p[\mu_0 + \epsilon \mu_1] = p_0 + \epsilon p_1 + O(\epsilon^2)$. Here, $\mu_0 \in D(E)$, and $\mu_1 \in T_{\mu} \mathscr{P}_2(X) \cap L^{\infty}(X)$ is an element of the tangent space, so that $\int_{X} \d \mu_1 = 0$.  Inserting this expansion into \eqref{eq:poisson1}, we obtain at zeroth-order
        \begin{align*}
            - \diver \left( \left(  \P[\mu_0] + r \mathbb{I}\right) \grad p_0 \right) = S,
        \end{align*}
        with zero-flux boundary conditions for $p_0$. Collecting terms of first order in $\eps$, we have
        \begin{align*}
            - \diver \left( \left( \P[\mu_0] + r \mathbb{I} \right) \grad p_1 + \P[\mu_1] \grad p_0 \right) = 0,
        \end{align*}
        again with zero-flux boundary conditions for $p_1$.
        Multiplication of this equation by $p_0$ and integration by parts gives
        \begin{align}
            \label{eq:identity1}
            \int_X \grad p_0 \cdot \C \grad p_1 \d \mu_0 + \int_{\Omega} r \grad p_0 \cdot \grad p_1 \d x = - \int_X \grad p_0 \cdot \C \grad p_0 \d \mu_1.
        \end{align}
Therefore, we have
        \begin{alignat*}{2}
            \leftindex_{T_{\mu}^* \mathscr{P}_2} {\left\langle \frac{\delta E[\mu_0]}{\delta \mu}, \mu_1 \right\rangle}_{T_{\mu} \mathscr{P}_2} &= &&\frac{\d}{\d\eps} E [\mu_0 + \eps \mu_1]\Bigr|_{\eps = 0} \notag \\ 
            & = &&\frac{\d}{\d\eps} \int_X \big[  \left( \grad p_0 + \eps \grad p_1 \right) \cdot \C \left( \grad p_0 + \eps \grad p_1 \right) + \frac{\nu}{\gamma} | \C |^{\gamma} \big] (\d \mu_0+ \eps \d \mu_1)\Bigr|_{\eps = 0} \notag \\
            & &&+\frac{\d}{\d\eps} \int_{\Omega}  r | \grad p_0 + \eps \grad p_1 |^2 \d x\Bigr|_{\eps = 0} \notag \\
            &= &&\int_X 2  \grad p_0 \cdot \C \grad p_1 \d \mu_0 + \int_X \left(  \grad p_0 \cdot \C \grad p_0 + \frac{\nu}{\gamma} | C |^{\gamma} \right) \d \mu_1 \notag \\
            & &&+ \int_{\Omega} 2  r \grad p_0 \cdot \grad p_1 \d x. \notag
        \end{alignat*}
        Using \eqref{eq:identity1} we finally obtain
        \begin{align} \label{eq:dEmu}
            \leftindex_{T_{\mu}^* \mathscr{P}_2} {\left\langle \frac{\delta E[\mu_0]}{\delta \mu}, \mu_1 \right\rangle}_{T_{\mu} \mathscr{P}_2} = \int_X \left( -  \grad p_0 \cdot \C \grad p_0 + \frac{\nu}{\gamma} | \C |^{\gamma} \right) \d \mu_1.
        \end{align}
    \end{proof}
\end{proposition}

We say that $\mu_0 = \mu_0 (x, \C)$ is a critical point of \eqref{eq:energy} if there exists $K \in \R$ such that the first variation \eqref{eq:first_variation} is constant, i.e., 
\begin{align}
    \label{eq:critical_point}
    - \grad p_0 \cdot \C \grad p_0 + \frac{\nu}{\gamma } | \C |^{\gamma} = K \quad \forall \, x \in \overline{\Omega}, \forall \, \C \in \overline{\mathscr{S}^d}(\R),
\end{align}
where $p_0$ solves \eqref{eq:poisson1}-\eqref{eq:bc_poisson} with the permeability $\P[\mu_0]$. Note also that if $\mu_0$ is a critical point of $E[\mu]$ and $\mu_1 \in T_{\mu} \mathscr{P}(X)$ then, by \eqref{eq:dEmu}, the function $e_{\mu_1}(\eps) : = E [\mu_0 + \eps \mu_1]$ satisfies $\left. \frac{\d e_{\mu_1}(\eps)}{\d\eps} \right|_{\eps=0}
 = 0$.

\begin{proposition}\label{no_critical_points}
    Let $S \not \equiv 0$. Then the energy \eqref{eq:energy} has no critical point in $X$.
    \begin{proof}
        For contradiction, suppose that $\mu_0$ is a critical point of $E[\mu]$. Then, there exists a constant $K \in \R$ such that \eqref{eq:critical_point} is satisfied. Now, we select $\C = \alpha \mathbb{I}$ for $\alpha \neq 0$. Thus,
        \begin{align*}
            | p_{x_i} |^2 = \frac{\nu \alpha^{\gamma-1}}{\gamma } - \frac{K}{\alpha} \quad \forall i = 1, \dots, d.
        \end{align*}
        Hence, $p_{x_i}$ is constant for all $i$, and an appropriate choice of the constant $\alpha\neq 0$ leads to a contradiction to the zero-flux boundary condition \eqref{eq:bc_poisson}.
    \end{proof}
\end{proposition}

Later we shall restrict $E$ to certain submanifolds of $\mathscr{P}_2(X)$ on which critical points exist.

Next, we establish strict convexity along straight lines
of the energy functional \eqref{eq:energy}, constrained by the Poisson equation \eqref{eq:poisson1}-\eqref{eq:bc_poisson}.

\begin{proposition}\label{prop:conv}
    The energy functional \eqref{eq:energy} constrained by \eqref{eq:poisson1}-\eqref{eq:bc_poisson} is strictly convex on its domain $D(E)$ and convex on $\mathscr{P}_2(X)$.
\end{proposition}
\begin{proof}
    We calculate the second-order variation of $E$ at $\mu_0 \in D(E)$. Let $\eps\in\R$ and $\mu_1 \in D(E) \cap T_{\mu} \mathscr{P}_2(X)$. We set $p[\mu_0 + \epsilon \mu_1] = p_0 + \epsilon p_1 + \frac{\epsilon^2}{2} p_2 + O(\epsilon^3)$. The Poisson equation \eqref{eq:poisson1} expanded at the zeroth, first, and second order in $\epsilon$ reads, respectively
    \begin{align*}
        - \diver \left( \left( \mathbb{P}[\mu_0] + r \mathbb{I} \right) \grad p_0 \right) &= S, \\
        - \diver \left( \left( \mathbb{P}[\mu_0] + r \mathbb{I} \right) \grad p_1 + \mathbb{P}[\mu_1] \grad p_0 \right) &= 0, \\
        - \diver \left( \frac{1}{2} \left( \mathbb{P}[\mu_0] + r \mathbb{I} \right) \grad p_2 + \mathbb{P}[\mu_1] \grad p_1 \right) &= 0,
    \end{align*}
    with no-flux boundary conditions \eqref{eq:bc_poisson} for $p_0$, $p_1$ and $p_2$.
    Multiplication of the zeroth order equation by $p_2$, the second order equation by $p_0$ and the first order equation by $p_1$, followed by integration by parts, gives
    \begin{align}
    \label{eq:comp_convex1}
        \int_{\Omega} S p_2 \d x &= \int_{\Omega} \grad p_2 \cdot \left( \mathbb{P}[\mu_0] + r \mathbb{I} \right) \grad p_0 \d x \notag \\
        &= - 2 \int_{\Omega} \grad p_0 \cdot \mathbb{P}[\mu_1] \grad p_1 \d x \notag \\
        &= 2 \int_{\Omega} \grad p_1 \cdot \left( \mathbb{P}[\mu_0] + r \mathbb{I} \right) \grad p_1 \d x.
    \end{align}
    Moreover, multiplication of \eqref{eq:poisson1} by $ p$ and integrating by parts yields
    \begin{align*}
        E[\mu] = \int_{\Omega}  S p[\mu] \d x + \int_X \frac{\nu}{\gamma} | \C |^{\gamma} \d \mu.
    \end{align*}
    Using the latter expression for the energy, we compute its second variation by employing the identity \eqref{eq:comp_convex1}. Therefore, we have
    \begin{align*}
        \leftindex_{T_{\mu}^* \mathscr{P}_2} {\left\langle \frac{\delta^2 E[\mu_0]}{\delta \mu^2} \mu_1, \mu_1\right\rangle}_{T_{\mu} \mathscr{P}_2} &= \frac{d^2}{d \eps^2} \left( \int_{\Omega}  S (p_0 + \eps p_1 + \frac{\eps^2}{2} p_2) \d x + \int_X \frac{\nu}{\gamma} | \C |^{\gamma} (\d \mu_0 + \eps \d \mu_1) \right) \Biggr|_{\eps = 0} \\
        &=  \int_{\Omega} S p_2 \d x \\
        &= 2  \int_{\Omega} \grad p_1 \cdot \left( \mathbb{P}[\mu_0] + r \mathbb{I} \right) \grad p_1 \d x \\
        &= 2  \int_X \grad p_1 \cdot \C \grad p_1 \d \mu_0 + 2  \int_{\Omega} r | \grad p_1 |^2 \d x > 0
    \end{align*}
    as $\C$ is positive semidefinite and $r > 0$. This, together with the convexity of the domain $D(E)$, concludes the proof.
\end{proof}

Contrary to the discrete model discussed in Section \ref{sec:discrete} and its related continuum limits, our model exhibits a distinctive feature - its convexity is not contingent upon the value of the metabolic exponent $\gamma>0$. Indeed, in the previously established models the energy functional is convex if and only if $\gamma\geq 1$ - see \cite{HV23} for the discrete case and \cite{kreusser2019, HMP15, HMPS16, CMS-2022, portaro2023} for the various versions of the continuum models.
In contrast, in our model the metabolic term in the energy \eqref{eq:energy} is merely the $\C$-moment of $\mu$ of order $\gamma>0$, and thus linear in $\mu$.
The nonlinearity of \eqref{eq:energy} comes from the kinetic energy term, which is convex according to Proposition \ref{prop:conv}.

We emphasize that the metabolic term $\int_{X} \frac{\nu}{\gamma} | \C |^{\gamma}  \d \mu$ in \eqref{eq:energy} exhibits geodesic convexity for $\gamma \ge 1$. Indeed, the integrand function is bounded below and $\lambda$-convex for $\gamma \ge 1$. Then, the geodesic convexity follows by \cite[Proposition 3.5]{ambrosio_handbook}. Thus, transitioning to the appropriate notion of convexity in our context, namely geodesic convexity, reiterates the essential condition $\gamma \ge 1$, as previously noted in the literature cited above.

Next we study the lower semicontinuity of the energy functional \eqref{eq:energy}.

\begin{theorem}\label{thm:lsc_D(E)}
    The energy functional $E$ defined in \eqref{eq:energy} is lower semicontinuous in $B:= \mathscr{P}_2(X) \cap L^2_{\mathrm{loc}} \left( \Omega ; L^1_{\mathrm{loc}} (\overline{\mathscr{S}^d}(\R)) \right)$.

    \begin{proof}
        Consider $\mu \in B$ and a sequence $\mu_n \in B$ converging to $\mu$ in the strong $L^2_{\mathrm{loc}} \left( \Omega ; L^1_{\mathrm{loc}} (\overline{\mathscr{S}^d}(\R)) \right)$ topology and in $\mathscr{P}_2(X)$. 
        Clearly the metabolic term $\int_X \frac{\nu}{\gamma} | \C |^{\gamma} \d \mu$ is lower semicontinuous and $\mathbb{P}[\mu_n]$ is absolutely continuous with respect to the Lebesgue measure on $\Omega$ for all $n$.
        We denote by $p_n$ the unique weak solution of the Poisson equation \eqref{eq:poisson1} with permeability tensor $\P[\mu_n]$ constructed in Lemma \ref{poisson_solvability}. The uniform bound \eqref{eq:estimate_gradp} implies that we can extract a subsequence converging weakly, say $p_n \weakconv p$ in $H^1(\Omega)$. Obviously, the background permeability term $\int_{\Omega}  r | \grad p |^2 \d x$ is lower semicontinuous.
        Moreover, multiplying the Poisson equation by $p_n = p[\mu_n]$, integrating by parts and using \eqref{eq:estimate_gradp}, we obtain the following estimate for the kinetic term
        \begin{align}
            \int_X  \grad p_n \cdot \C \grad p_n \d \mu_n &= \int_X  \big{|} \sqrt{\C} \grad p_n \big{|}^2 \d \mu_n \notag \notag \\
            &= \int_{\Omega} S p_n \d x - r \int_{\Omega} |\grad p_n|^2 \d x \notag \\
            &\le \frac{C_{\Omega}^2}{r} \| S \|_{L^2(\Omega)}.
            \label{eq:ubdd}
        \end{align}
        Then, \cite[Theorem 4.6]{ambrosio_handbook} ensures the existence of limit points in the narrow topology of $\sqrt{\C} \grad p_n \mu_n$ as $n \rightarrow \infty$ and we identify the limit as $\sqrt{\C} \grad p \mu$ due to the weak convergence $\grad p_n \weakconv \grad p$ in $L^2(\Omega)$ and the strong convergence $\mu_n \rightarrow \mu$ in $L^2_{\mathrm{loc}} \left( \Omega ; L^1_{\mathrm{loc}} (\overline{\mathscr{S}^d}(\R)) \right)$. Moreover, the Theorem gives the lower semincontinuity 
        \begin{align*}
            \int_X  \big{|} \sqrt{\C} \grad p \big{|}^2 \d \mu \le \liminf_{n \rightarrow \infty} \int_X  \big{|} \sqrt{\C} \grad p_n \big{|}^2 \d \mu_n.
        \end{align*}
        We conclude
        \begin{align}
            \label{eq:proof_bound}
            \int_{X} \left( \grad p \cdot \C \grad p + \frac{\nu}{\gamma} | \C |^{\gamma} \right) \d \mu + r \int_{\Omega} | \grad p |^2 \d x \le \liminf_{n \rightarrow \infty} E[\mu_n].
        \end{align}
        Moreover,
        \cite[Theorem 4.6]{ambrosio_handbook} gives $\int_X \xi (x, \C) \cdot \sqrt{C} \grad p_n \d \mu_n \rightarrow \int_X \xi(x, \C) \cdot \sqrt{C} \grad p \d \mu$ for every continuous vector test function $\xi$ with at most linear growth in $\C$. Choosing $\xi(x, \C) = \varphi(x) \sqrt{\C}$ with $\varphi \in C^{\infty}(\overline{\Omega})$ shows
        \begin{align*}
            \int_X \varphi(x) \cdot \C \grad p_n \d \mu_n = \int_{\Omega} \varphi(x) \cdot \mathbb{P}[\mu_n] \grad p_n \d x
        \end{align*}
        and passing to the limit, we have
        \begin{align*}
            \int_X \varphi(x) \cdot \C \grad p \d \mu = \int_{\Omega} \varphi(x) \cdot \P[\mu] \grad p \d x.
        \end{align*}
        This allows us to pass to the limit in the weak formulation of the Poisson equation and we conclude that $p$ is the unique solution in $H$ of \eqref{eq:poisson1}-\eqref{eq:bc_poisson}. 
        Lower semicontinuity of $E$ in $B$ with respect to its strong $L^2_{\mathrm{loc}} \left( \Omega ; L^1_{\mathrm{loc}} (\overline{\mathscr{S}^d}(\R)) \right)$ topology follows.
    \end{proof}
\end{theorem}
One can immediately notice though that lower semicontinuity fails on the complement of the set $\tilde{D}(E) := \{ \mu \in \mathscr{P}_2(X) : \P[\mu] \in L^1(\Omega) \; \mathrm{or} \; \int_{X} | \C |^{\gamma} \d \mu = \infty \}$. Indeed, suppose $\mu$ is chosen such that $\P[\mu]$ is not absolutely continuous with respect to the Lebesgue measure and $\int_{X} | \C |^{\gamma} \d \mu < \infty$. Then, $E[\mu] = \infty$, and there exists a sequence $\mu_n\in\mathscr{P}_2(X)$ with $\P[\mu_n]$ absolutely continuous with respect to the Lebesgue measure, and such that $W_2(\mu_n, \mu) \rightarrow 0$. Now, assuming lower semicontinuity $E[\mu] \le \liminf_{n \rightarrow \infty} E[\mu_n]$ immediately leads to a contradiction because $E[\mu]=+\infty$ and
\begin{align*}
    \liminf_{n \rightarrow \infty} E[\mu_n] < \infty
\end{align*}
since the kinetic term in the energy functional is uniformly bounded by \eqref{eq:ubdd}.
This situation is typical for functionals that do not grow superlinearly, as discussed in \cite{santambrogio2015} (see section 7.1.2). In particular our energy grows linearly only due to the metabolic term.

We note that lower semicontinuity is a typical condition for establishing well-posedness of Wasserstein gradient flows via the minimizing movements scheme by De Giorgi in \cite{degiorgi}. We shall investigate this problem in an upcoming work. 

To establish global lower semicontinuity, it is common to compute a relaxation $\overline{E} : \mathscr{P}_2(X) \rightarrow \R \cup \{ \infty \}$ of the functional, defined as
\begin{align}
    \label{eq:def_energy_relaxed}
    \overline{E}[\mu] := \sup \left\{ G : \mathscr{P}_2(X) \rightarrow \R \cup \{ \infty \}; \; G \le E \, \, \mathrm{and} \, \, G \, \, \mathrm{is} \, \, \mathrm{narrowly} \, \,\mathrm{lower} \, \, \mathrm{semicontinuous}\right\}.
\end{align}
Then $\overline{E}$ is lower semicontinous since it is the supremum of a family of lower semicontinuous functionals. Note that $\overline{E}$ coincides with $E$ within the domain $\tilde{D}(E)$. We also have the representation formula
\begin{align}
    \label{eq:energy_relaxed}
    \overline{E}[\mu] = \inf \left\{ \liminf_{n \to \infty} E[\mu_n]; \;  \mu_n \rightarrow \mu \, \, \mathrm{narrowly} \right\}.
\end{align}
\begin{proposition}
The relaxed functional \eqref{eq:energy_relaxed} assumes its minimum as the infimum of the functional \eqref{eq:energy}, that is
\begin{align*}
    \min_{\mu \in \mathscr{P}_2(X)} \overline{E}[\mu] = \inf_{\mu \in \mathscr{P}_2(X)} E[\mu].
\end{align*}
\begin{proof}
    We proceed along the lines of the proof of \cite[Theorem 7.5]{rindler}. Take an infimizing sequence $\{ \mu_k \} \subset \mathscr{P}_2(X)$ such that
    \begin{align*}
        \lim_{k \rightarrow \infty} \overline{E}[\mu_k] = l := \inf_{\mu} \overline{E}[\mu] < \infty.
    \end{align*}
    Then, there exists a constant $C >0$ such that $\overline{E}[\mu_k] \le C$ for all $k \in \N$. From this bound, we find that the $|\C|^{\gamma}$ moment of $\mu_k$ is uniformly bounded in $k \in \N$, implying that the sequence $\mu_k$ is tight and therefore (by Prokhorov's Theorem \cite{prokhorov}) narrowly compact.
    Hence, we can select a subsequence, again denoted by $\mu_k$, such that $\mu_k \rightarrow \mu_*$ narrowly in $\mathscr{P}(X)$ for some $\mu_* \in \mathscr{P_2}(X)$. Then, by lower semicontinuity of $\overline{E}$ in the narrow topology in $\mathscr{P}(X)$, we have
    \begin{align*}
        l \le \overline{E}[\mu_*] \le \liminf_{k \rightarrow \infty} \overline{E}[\mu_k] = l.
    \end{align*}
    Consequently, $\overline{E}[\mu_*] = l$ and $\mu_*$ is the minimizer.
    Finally, to conclude, we notice that $\inf_{\mu} E \le \overline{E} \le E$ from the definition of the relaxed functional \eqref{eq:def_energy_relaxed}.
\end{proof}
\end{proposition}

\begin{remark}
If $S \not \equiv 0$, then
\begin{align*}
    \inf_{\mu \in \mathscr{P}_2(X)} E[\mu] > 0.
\end{align*}
Indeed, assume that $\inf_{\mu \in \mathscr{P}_2(X)} E[\mu] = 0$ and take an infimizing sequence $\{ \mu_k \} \subset \mathscr{P}_2(X)$ such that $\P[\mu_k] \in L^1(\Omega)$ for all $k \in \N$. Then $\lim_{k \rightarrow \infty} E[\mu_k] = 0$ implies that $\int_{\Omega}  r | \grad p_k |^2 \rightarrow 0$ as $k\to\infty$. This in turn means that $\lim_{k \rightarrow \infty} p_k$ is a constant, which violates the boundary condition \eqref{eq:bc_poisson}.
\end{remark}

\subsection{Wasserstein gradient flows} \label{wasserstein}
We calculate the $\mathscr{P}_2(X)$-Wasserstein gradient flow of the energy functional \eqref{eq:energy} constrained by \eqref{eq:poisson1}--\eqref{eq:bc_poisson}. The resulting system of partial differential equations can be seen as a model that accounts for dynamic changes in the transportation network.


\begin{proposition}\label{wass_GF}
    The formal $\mathscr{P}_2(X)$-Wasserstein gradient flow of the energy functional \eqref{eq:energy} constrained by \eqref{eq:poisson1}-\eqref{eq:bc_poisson} is given by
    \begin{align}
        \label{eq:wass_mu}
        \partial_t \mu + \diver_x \left( 2 \mu D^2 p \C \grad p \right) + \diver_{\C} \left( \mu \left( \grad p \otimes \grad p - \nu | \C |^{\gamma-2} \C \right) \right) = 0,
    \end{align}
    for $x\in\Omega$, $\C\in\mathscr{S}^d(\R)$ and $t\geq 0$,
    with the boundary conditions
    \begin{align*}
        \mu \grad_x \left( \grad p \cdot \C \grad p \right) \cdot n = 0 \quad &\mathrm{on} \, \, \partial \Omega \times \mathscr{S}^d(\R) \times (0, \infty), \\
         \mu \mathsf{v}\left(  \grad p \otimes \grad p - \nu | \C |^{\gamma - 2} \C \right) \cdot n_{\mathscr{S}} = 0 \quad &\mathrm{on} \, \, \Omega\times\partial \mathscr{S}^d(\R) \times (0, \infty),
    \end{align*}
    where $n_{\mathscr{S}}$ is the unit outer normal vector to $\partial \mathscr{S}^d(\R)$ and $n$ is the unit outer normal vector to $\partial \Omega$.
    \begin{proof}
    
    The formal Wasserstein gradient flow is given by \eqref{eq:gradient_structure} where both the divergence and gradient operators are applied in the phase space $X$, i.e., with respect to $x$ and the tensor $\C$. We have
    \begin{align*}
        \left( \grad_{\C} \left( - \grad p \cdot \C \grad p + \frac{\nu}{\gamma} | \C |^{\gamma} \right) \right)_{kl} &= \frac{\partial}{\partial C_{kl}} \left( - p_{x_i} C_{ij} p_{x_j} + \frac{\nu}{\gamma} | C_{ij} |^{\gamma} \right) \notag \\
        &= - p_{x_k} p_{x_l} + \nu | C_{kl} |^{\gamma - 2} C_{kl} \notag \\
        &= \left( - \grad p \otimes \grad p + \nu | \C |^{\gamma-2} \C \right)_{kl} \quad \forall \, k, l = 1, \dots, d,
    \end{align*}
    where we have used the Einstein summation notation with the indices $i$, $j$. Moreover, we evaluate
    \begin{align*}
        \grad_{x} \left( - \grad p \cdot \C \grad p + \frac{\nu}{\gamma} | \C |^{\gamma} \right) = - 2  D^2 p \C \grad p,
    \end{align*}
    where $D^2 p$ denotes the Hessian of $p$. Consequently, inserting those terms into \eqref{eq:gradient_structure} we obtain \eqref{eq:wass_mu} with no flux boundary conditions \cite[Remark 6.1]{ambrosio_handbook}.
    \end{proof}
\end{proposition}
Applying \eqref{eq:flow_dissipation} to the $2$-Wasserstein flow, we recollect the energy dissipation law
\begin{align*}
    \frac{d}{d t}E[\mu(t)] = - \int_{X} \left( 4 | D^2 p \C \grad p |^2 + | \grad p \otimes \grad p - \nu | \C |^{\gamma-2} \C |^2 \right) \d \mu \le 0.
\end{align*}

Interpreting \eqref{eq:wass_mu} as a hyperbolic transport equation for the probability measure $\mu = \mu(x, \C)$, we have the following characteristic system for x and \C,
    \begin{align}
        \label{eq:wass_xt}
        \frac{dx}{dt} &=  \grad_x \left( \grad p \cdot \C \grad p \right), \\
        \label{eq:wass_Ct}
        \frac{d \C}{dt} &=  \grad p \otimes \grad p - \nu | \C |^{\gamma - 2} \C,
    \end{align}
which is coupled to the equation for $\mu$ along the characteristics,
\begin{equation}
    \label{eq:wass_mu2}
    \frac{\d \mu}{\d t} = - \mu \bigg[ \Delta_x \left( \grad p \cdot \C \grad p \right) - \nu \left( \gamma - 1 \right) | \C |^{\gamma-2} \bigg].
\end{equation}
Let $\mu_0 \in \mathscr{P}_2(X)$, if there is sufficient regularity of the characteristic flow map $\mathscr{X}$ on $X \times [0, T]$, then $\mu := \left( \mathcal{X} \right)_{\#} \mu_0 $ is a solution of \eqref{eq:wass_mu2} in $[0, T]$ \cite[Lemma 8.1.6]{ambrosiogigli}.

For a monokinetic $\mu$, i.e., $\mu(x, \C) = \rho(x) \delta \left( \C - \hat{\C}(x) \right)$, with $\rho \in \mathscr{P}(\overline{\Omega})$, the energy takes the form
\begin{align}
    \label{eq:monokinetic_energy}
    E[\rho, \hat{\C}] = \int_{\Omega} \bigg[ \left(  \grad p \cdot \hat{\C} \grad p + \frac{\nu}{\gamma} | \hat{\C} |^{\gamma} \right) + r | \grad p|^2 \bigg] \d \rho.
\end{align}
This is a generalization of the energy functional considered in \cite{CMS-2022}. We now compute the $\mathscr{P}_2(X)$-Wasserstein gradient flow with respect to $\rho$ and the $L^2$ flow with respect to $\hat{\C}$ of the above energy functional.
\begin{proposition}
    The formal Wasserstein gradient flow in $\rho$ and $L^2$-gradient flow in $\hat{\C}$ of the energy functional \eqref{eq:monokinetic_energy} constrained by
    \begin{align}
        \label{eq:monokinetic_poisson}
        - \diver \left( \left( \rho \hat{\C} + r \mathbb{I} \right) \grad p \right) = S, \quad &\mathrm{in} \, \, \Omega \\ 
        n \cdot \left( \rho \hat{\C} + r \mathbb{I} \right) \grad p = 0 \quad &\mathrm{on} \, \, \partial \Omega \notag
    \end{align}
    is given by
    \begin{align}
        \label{eq:monokinetic_wass_rho}
        \partial_t \rho + \diver \bigg[ \rho \grad \left( \grad p \cdot \hat{\C} \grad p - \frac{\nu}{\gamma} | \hat{\C} |^{\gamma} \right) \bigg] = 0 \quad &\mathrm{in} \, \, \Omega \times (0, \infty), \\
        \label{eq:monokinetic_L2_C}
        \partial_t \hat{\C} = \rho \left( \grad p \otimes \grad p - \nu | \hat{\C} |^{\gamma -2} \hat{\C} \right) \quad &\mathrm{in} \, \, \Omega \times (0, \infty)
    \end{align}
    with the boundary condition
    \begin{align*}
        \rho \grad \left( \grad p \cdot \hat{\C} \grad p - \frac{\nu}{\gamma} | \hat{\C} |^{\gamma} \right) \cdot n = 0 \quad &\mathrm{on} \, \, \partial \Omega \times (0, \infty). \\
    \end{align*}
    \begin{proof}
    The conservation law \eqref{eq:monokinetic_wass_rho} is derived in analogy to what has been done in Proposition \ref{wass_GF}, where the gradient and divergence operators are applied solely with respect to $x \in \Omega$.

    To compute the $L^2$-gradient flow with respect to $\hat{\C}$, we expand $\hat{\C} = \hat{\C}_0 + \eps \hat{\C}_1$ and $p [\hat{\C}] = p_0 + \eps p_1 + O(\eps^2)$ with $\eps \in \R$, where $\hat{\C}_0$ is a symmetric positive definite tensor and $\hat{\C}_1$ is symmetric.
    Inserting this expansion into \eqref{eq:monokinetic_poisson}, we obtain at zeroth-order
    \begin{align*}
        - \diver \left( \left( \rho \hat{\C}_0 + r \mathbb{I} \right) \grad p_0 \right) = S,
    \end{align*}
    with zero-flux boundary conditions for $p_0$. On the other hand, collecting terms of first order in $\eps$, we have
    \begin{align*}
        -\diver \left( \left( \rho \hat{\C}_0 + r \mathbb{I} \right) \grad p_1 + \rho \hat{\C}_1 \grad p_0 \right) = 0. 
    \end{align*}
    Multiplying the zeroth-order Poisson equation by $p_1$, the first-order one by $p_0$ and integrating by parts we get
    \begin{align}
    \label{eq:monokinetic_identity}
        \int_{\Omega} S p_1 \d x = \int_{\Omega} \left( \rho \grad p_1 \cdot \hat{\C}_0 \grad p_0 + r \grad p_1 \cdot \grad p_0 \right) \d x = - \int_{\Omega} \rho \grad p_0 \cdot \hat{\C}_1 \grad p_0 \d x.
    \end{align}
    Note that, a multiplication of \eqref{eq:monokinetic_poisson} by p and subsequent integration by parts leads to the equivalent form of the energy \eqref{eq:monokinetic_energy}
    \begin{align*}
        E[\rho, \hat{\C}] = \int_{\Omega} \left( S p[\rho, \hat{\C}] + \frac{\nu}{\gamma} | \hat{\C} |^{\gamma} \right) \d x.
    \end{align*}
    Using the latter expression, we compute its first variation in $\hat{\C}$ by employing identity \eqref{eq:monokinetic_identity}. Thus,
    \begin{align*}
         \leftindex_{L^2(\Omega; S_d)} {\left( \frac{\delta \E[\rho, \hat{\C}_0]}{\delta \hat{\C}}, \hat{\C}_1 \right)}_{L^2(\Omega; S_d)} &= \frac{\d}{\d \eps} E[\rho, \hat{\C}_0 + \eps \hat{\C}_1]\Bigr|_{\eps = 0} \\
        &= \int_{\Omega} \left( S p_1 + \nu \rho | \hat{\C}_0 |^{\gamma-2} \hat{\C}_1 : \hat{\C}_0 \right) \d x \\
        &= - \int_{\Omega} \rho \hat{\C}_1 : \left( \grad p_0 \otimes \grad p_0 - \nu | \hat{\C}_0 |^{\gamma-2} \hat{\C}_0 \right) \d x,
    \end{align*}
    where $S_d$ is the space of symmetric real $d \times d$ matrices. Therefore, \eqref{eq:monokinetic_L2_C} is proved.
    \end{proof}
\end{proposition}
We remark that \eqref{eq:monokinetic_L2_C} with uniform measure $\rho$ is precisely the model considered in \cite{CMS-2022}.


\subsection{Reduced Wasserstein gradient flow}\label{subsec:reduced}

Next, we compute the reduced Wasserstein gradient flow of the energy functional \eqref{eq:energy}, constrained by \eqref{eq:poisson1}--\eqref{eq:bc_poisson}, with respect to $\C\in \overline{\mathscr{S}^d}(\R)$. This corresponds to a change of paradigm, since $\mu$ will now be considered as a map with domain $\overline{\Omega}$ into the $2$-Wasserstein space on positive semidefinite matrices. This choice is motivated by the discrete model outlined in Section \ref{sec:discrete}. In particular, in the discrete model the energy functional \eqref{eq:discrete_energy} is defined in terms of edge conductivities on a graph fixed in the physical space.

It is a priori unclear if this reduced $\C$-flow constitutes a formal energy-dissipating gradient flow. However, in the subsequent discussion we shall demonstrate that this is the case.
We start by defining the new gradient structure triple $\left( N, E, \mathbb{H} \right)$, where $N := C^{\infty}_{+} \left( \overline{\Omega}; \mathrm{Prob}_{+,2} (\overline{\mathscr{S}^d}(\R) \right)$. Here, $\mathrm{Prob}_{+,2}$ denotes the space of all smooth, positive probability densities with finite second moments.
Obviously $T_{\mu} N$ can be identified with  $ \big{\{} s : X \rightarrow \R \; \mathrm{smooth} : \int_{\overline{\mathscr{S}^d}(\R)} s(x, \C) \d \C = 0 \; \forall x \in \overline{\Omega} \big{\}}$. Note that the position density $\mu_{\C} (x) = 1$ in $\Omega$ for $\mu \in N$.
We now define the Riemannian structure $\mathbb{H}(\mu) : T_{\mu} N \rightarrow T_{\mu}^* N$ by $\mathbb{H}(\mu) s = \Phi$ if and only if $\Phi$ is the (weak) solution of
\begin{align*}
    \begin{cases}
        - \diver_{\C} \left( \mu(x, \C) \grad_{\C} \Phi(x, \C) \right) = s(x, \C) \quad &\mathrm{in} \; \mathscr{S}^d(\R) \\
        \mu \grad_{\C} \Phi \cdot n_{\mathscr{S}^d(\R)} = 0 \quad &\mathrm{on} \; \partial \mathscr{S}^d(\R)
    \end{cases}
\end{align*}
for all $x \in \overline{\Omega}$.
We compute the inverse $\mathbb{L}(\mu) := \mathbb{H}(\mu)^{-1} : D(\mathbb{L}(\mu)) \subseteq T_{\mu}^* N \rightarrow T_{\mu} N$ as $\mathbb{L}(\mu) \Phi = -\diver_{\C} \left( \mu \grad_{\C} \Phi \right)$ and identify $D(\mathbb{L}(\mu)) = \{ \Phi \in T^*_{\mu} N : \mu \grad_{\C} \Phi \cdot n_{\mathscr{S}^d(\R)} = 0 \; \mathrm{on} \; \partial \mathscr{S}^d(\R) \; \forall x \in \overline{\Omega} \}$ where $T^*_{\mu} N = \{ \Phi:X\rightarrow \R \; \mathrm{smooth} \} / C^{\infty}(\overline{\Omega})$.
Next, let us consider $s_1, s_2 \in T_{\mu} N$ and compute $\Phi_1$, $\Phi_2$ from 
\begin{align*}
    - \diver_{\C} \left ( \mu \grad_{\C} \Phi_1 \right) = s_1, \\
    - \diver_{\C} \left ( \mu \grad_{\C} \Phi_2 \right) = s_2,
\end{align*}
subject to no-flux boundary conditions. Then, the metric tensor $h_{\mu}(s_1, s_2) := \leftindex_{{T^*_\mu} N} {\left\langle \mathbb{H}(\mu) s_1, s_2 \right\rangle}_{{T_\mu} N}$ induced by $\mathbb{H}(\mu)$ is written as
\begin{align}
    h_{\mu}(s_1, s_2) &= \leftindex_{T^*_{\mu} N} {\left\langle \mathbb{H}(\mu) s_1, s_2 \right\rangle}_{T_{\mu} N} \notag \\
    &= \leftindex_{T^*_{\mu} N} {\left\langle \Phi_1, s_2\right\rangle}_{T_{\mu} N} \notag \\
    &= \int_{X} \Phi_1 \left( - \diver_{\C} \left( \mu \grad_{\C} \Phi_2 \right) \right) \d x \d \C \notag \\
    &= \int_X \grad_{\C} \Phi_1 \cdot \grad_{\C} \Phi_2 \d \mu,
    \label{eq:metric_tensor}
\end{align}
which is obviously a positive and symmetric tensor.

Now, let $t \in [0, 1]$ and let $\Phi = \Phi(x, \C, t)$ be a parametrized solution of the following system
\begin{eqnarray}
    \label{eq:parametrization}
        \left\{ 
        \begin{array}{ll} \displaystyle
        \partial_t \mu + \diver_{\C} \left( \mu \grad_{\C} \Phi \right) = 0 \\ \displaystyle
        \mu(x, \C, 0) = \mu_{0} (x, \C) \\ \displaystyle
        \mu(x, \C, 1) = \mu_{1} (x, \C),
        \end{array}
        \right.
    \end{eqnarray}
with zero-flux boundary condition. Here, $\mu(x, \C, t)$ represents a curve in $N$ connecting $\mu_0$ and $\mu_1$. We conclude that the distance on $N$ is given by 
\begin{align}
    \label{eq:distance}
    d_N  \left( \mu_0, \mu_1 \right) := \left( \int_{\Omega} W_2^2 \left( \mu_0 (x, \cdot), \mu_1 (x, \cdot) \right) \d x \right)^{\frac{1}{2}}
\end{align}
has the metric tensor $h_{\mu}$. Actually
\begin{align*}
    W_2^2 \left( \mu_0 (x, \cdot), \mu_1 (x, \cdot) \right) &= \min_{\Phi(x, \cdot)} \int_0^1 \int_{\overline{\mathscr{S}^d}(\R)} | \grad_{\C} \Phi \left( x, \C, t \right) |^2 \mu(x, \d \C, t) \d t \\
    &\le \int_0^1 \int_{\overline{\mathscr{S}^d}(\R)} | \grad_{\C} \Phi \left( x, \C, t \right) |^2 \mu(x, \d \C, t) \d t, \\
\end{align*}
for all $\Phi(x, \cdot)$ satisfying \eqref{eq:parametrization} where the first identity is the Benamou-Brenier Formula \cite{benamou-brenier}. Consequently, integrating over $\Omega$, we find
\begin{align*}
    d_N^2 \left( \mu_0, \mu_1 \right) &= \min_{\Phi} \int_0^1 \int_X | \grad_{\C} \Phi \left( x, \C, t \right) |^2 \mu(\d x, \d \C, t) \d t \\
    &= \min_{
    \mu(0) = \mu_0, \mu(1)=\mu_1} \int_0^1 h_{\mu(t)} \left( \partial_t \mu, \partial_t \mu \right) \d t.
\end{align*}
The identity follows since the square of the Riemannian distance of two points, i.e., $d^2_N$, is the minimum of the energy on the set of minimal length geodesics connecting those two points \cite[Section 9.2]{docarmo}. We conclude that $h_{\mu}$ is the metric tensor of the distance \eqref{eq:distance}.

We define the space $Q := L^2 \left( \Omega; \mathscr{P}_2 \left( \overline{\mathscr{S}^d}(\R) \right) \right)$ which we endow with the metric $d_N$ \eqref{eq:distance}. We remark that $\mu \in Q$ if and only if $\mu(x, \cdot) = \mu_x$ is a probability measure in $\mathscr{P}_2 \left( \overline{\mathscr{S}^d} (\R) \right)$ for Lebesgue a.e. $x \in \Omega$ such that $\int_{\Omega} \left( \int_{\overline{\mathscr{S}^d}(\R)} | \C |^2 \mu_x (\d \C) \right) \d x < \infty$ and $x \rightarrow \mu_x(B)$ is a Borel measurable function on $\overline{\Omega}$ for every Borel subset $B \subseteq \overline{\mathscr{S}^d}(\R)$.
We now redefine the functional as $E : Q \rightarrow [0, \infty]$ with $D_Q(E) := \big{\{} \mu \in Q : \int_X |\C|^{\gamma} \d \mu < \infty \big{\}}$. We take $\mu_0, \mu_1 \in Q$ and apply the Kantorovich-Rubinstein theorem \cite[Theorem 1.14]{villani2021}
\begin{align*}
    W_1\left( \mu_0(x, \cdot), \mu_1(x, \cdot) \right) = \sup_{\varphi=\varphi\left(\C\right)} \bigg{\{} &\int_{\overline{\mathscr{S}^d}(\R)} \varphi \left( \mu_0 (x, \d \C) - \mu_1 (x, \d \C) \right) : \; \| \varphi \|_{\mathrm{Lip}} \le 1 \bigg{\}}.
\end{align*}

Thus, since $\mathrm{Lip} (\varphi) = 1$ for $\varphi(\C) = \C$,
\begin{align*}
    \bigg{|} \int_{\overline{\mathscr{S}^d}(\R)} \C \left( \mu_0 (x, \d \C) - \mu_1 (x, \d \C) \right) \bigg{|} = \big{|} \mathbb{P}[\mu_0](x) - \mathbb{P}[\mu_1](x) \big{|} &\le W_1 \left( \mu_0(x, \cdot), \mu_1(x, \cdot) \right) \\ 
    &\le W_2 \left( \mu_0(x, \cdot), \mu_1(x, \cdot) \right)
\end{align*}
and
\begin{align}
    \label{eq:PL2_ineq}
    \| \mathbb{P}[\mu_0] - \mathbb{P}[\mu_1] \|_{L^2(\Omega)} \le d_N (\mu_0, \mu_1).
\end{align}
Therefore, $\mathbb{P}[\mu] \in L^2(\Omega) \subset L^1(\Omega)$ for all $\mu \in Q$ and $E[\mu] < \infty$ if and only if $\mu \in D_Q(E)$. 

We conclude:
\begin{proposition}
    The energy functional $E : D_Q(E) \rightarrow \R^+ $ defined in \eqref{eq:energy} is lower semicontinuous in $Q$.
    \begin{proof}
        Since the metabolic term is lower semicontinuous in the weak-* topology on $\mathscr{M}^+(X)$ it suffices to choose $\mu_n, \mu \in D_Q(E)$ with $d_N (\mu_n, \mu) \rightarrow 0$. We denote by $p_n := p[\mu_n]$ the unique weak solution of the Poisson equation \eqref{eq:poisson1} with permeability tensor $\mathbb{P}[\mu_n]$ constructed in Lemma \ref{poisson_solvability}. From \eqref{eq:PL2_ineq} we conclude that $\mathbb{P}[\mu_n] \rightarrow \mathbb{P}[\mu]$ in $L^2(\Omega)^{d \times d}$ strongly and by passing to the limit $n \rightarrow \infty$ in the Poisson equation we get $p := p[\mu] = \lim_{n \rightarrow \infty} p_n$. 

        We obtain, using Lemma $3.2$ and Proposition $3.3$ of \cite{CMS-2022}
        \begin{align*}
            \int_{\Omega} \grad p \cdot \mathbb{P}[\mu] \grad p \d x \le \liminf_{n \rightarrow \infty} \int_{\Omega} \big{|} \sqrt{\mathbb{P}[\mu_n]} \grad p_n \big{|}^2 \d x = \liminf_{n \rightarrow \infty} \int_{\Omega} \grad p_n \cdot \mathbb{P}[\mu_n] \grad p_n \d x,
        \end{align*}
        which yields the lower semicontinuity of $E$.
    \end{proof}
\end{proposition}

We arrive at the following representation of the reduced Wasserstein gradient flow.

\begin{proposition} \label{wasstype_GF}
    The reduced Wasserstein gradient flow of the energy functional \eqref{eq:energy} constrained by \eqref{eq:poisson1}-\eqref{eq:bc_poisson} is given by
    \begin{align}
        \label{eq:wasstype_mu}
        \partial_t \mu + \diver_{\C} \left( \mu \left(  \grad p \otimes \grad p - \nu | \C |^{\gamma-2} \C \right) \right) = 0 \quad \mathrm{in} \, \, \left( \overline{\Omega} \times \overline{\mathscr{S}^d}(\R) \right) \times (0, \infty),
    \end{align}
    with the boundary condition
    \begin{align*}
        \mu \mathsf{v} \left(  \grad p \otimes \grad p - \nu | \C |^{\gamma - 2} \C \right) \cdot n_{\mathscr{S}} = 0 \quad \mathrm{on} \, \, \partial \mathscr{S}^d(\R) \times (0, \infty).
    \end{align*}
    \end{proposition}
    Note that the initial datum $\mu$ for \eqref{eq:wasstype_mu} has to satisfy $\int_{\overline{\mathscr{S}^d}(\R)} \mu(x, \d \C) = 1$ for a.e. $x \in \Omega$. At the end of this Section we shall generalize this assumption.
    \begin{proof}
        The proof follows the one of Proposition \ref{wass_GF} where the divergence and gradient operator are applied solely with respect to $\C\in\overline{\mathscr{S}^d}(\R)$.
    \end{proof}

Applying \eqref{eq:flow_dissipation} to the reduced Wasserstein flow, we recollect the energy dissipation law
\begin{align*}
    \frac{\d}{\d t}E[\mu(t)] = - \int_{X} \big| \left(  \grad p \otimes \grad p - \nu | \C |^{\gamma-2} \C \right) \big|^2 \d \mu \le 0.
\end{align*}

Interpreting \eqref{eq:wasstype_mu} as a hyperbolic transport equation for the probability measure $\mu = \mu(x, \C)$, we have the following characteristic system,
    \begin{align}
        \label{eq:wasstype_Ct}
        \frac{d\C}{dt} &=  \grad p \otimes \grad p - \nu | \C |^{\gamma - 2} \C, \\
        \frac{dx}{dt} &= 0, \notag
    \end{align}
which is coupled to the equation for $\mu$ along the characteristics,
\begin{align}
    \label{eq:wasstype_mu2}
    \frac{\d \mu}{\d t} = \nu ( \gamma - 1) \mu | \C |^{\gamma-2}.
\end{align}
We notice that the measure $\mu$ decays along characteristics for $\gamma < 1$, is constant for $\gamma = 1$ and grows along characteristics for $\gamma > 1$.
Again, take $\mu_0 \in \mathscr{P}_2(X)$, if there is sufficient regularity on the characteristic flow map $\mathscr{C}$, then $\mu := \left( \mathscr{C} \right)_{\#} \mu_0 $ is a solution of \eqref{eq:wasstype_mu2} \cite[Lemma 8.1.6]{ambrosiogigli}. 

\begin{remark} \label{remark2}
    We can also revert back to the original set-up of probability measures on $X$ by redoing the formal Riemannian calculus of the beginning of Section \ref{subsec:reduced} replacing the space $N$ by $M$ of Section \ref{subsec:GFWass}. Then \eqref{eq:wasstype_mu} is still the formal gradient flow equation, where the solution $\mu = \mu(t)$ is now sought in $\mathscr{P}_2(X)$, for given initial datum in $\mathscr{P}_2(X)$, instead of Q.

    In order to better understand the dynamics of the flow we start by foliating the space $\mathscr{P}_2(X)$ according to the position density measure $\mu_{\C}$. We define the leaves of the foliation, for given $\rho \in \mathscr{P}(\overline{\Omega})$
    \begin{align*}
        Q_{\rho} = \{ \mu \in \mathscr{P}_2(X) : \mu_{\C} = \rho \}.
    \end{align*}
    Note that by the disintegration theorem of probability measures \cite[Theorem 5.3.1]{ambrosiogigli} we can write
    \begin{align*}
        \mu = \int_{\overline{\Omega}} \mu_x (\C) \rho (\d x),
    \end{align*}
    where $\mu_x$ is a probability measure in $\mathscr{P}\left(\overline{\mathscr{S}^d}(\R)\right)$ defined uniquely $\rho$-a.e. on $\overline{\Omega}$ such that the map $x \rightarrow \mu_x (B)$ is Borel-measurable on $\overline{\Omega}$ for every Borel-measurable subset $B \subseteq \overline{\mathscr{S}^d}(\R)$ and $\int_{\overline{\Omega}} \int_{\overline{\mathscr{S}^d}(\R)} | \C |^2 \mu_x(\d \C) \d \rho < \infty$. Thus, we can identify $Q_{\rho}$ with the space of all these probability measures $\mu_x(d \C)$. Note that the reduced Wasserstein flow of Proposition \ref{wasstype_GF} leaves $Q_{\rho}$ invariant for every $\rho \in \mathscr{P}(\overline{\Omega})$.

    By proceeding as above (where $\mathrm{vol}(\Omega) \rho$ is the Lebesgue measure on $\Omega$) we conclude that the metric on $Q_{\rho}$ associated with the (gradient) flow of Proposition \ref{wasstype_GF} is given by
    \begin{align*}
        d_{Q_{\rho}} (\mu_0, \mu_1 ) = \left( \int_{\Omega} W_2^2 (\mu_{0,x}, \mu_{1,x}) \d \rho \right)^{\frac{1}{2}}.
    \end{align*}
    Here $\mu_{0,x}, \mu_{1,x}$ are the disintegrations of $\mu_0$ and $\mu_1$, respectively, with respect to the position density measure $\rho$. From the disintegration theorem we conclude
    \begin{align*}
        \mathbb{P}[\mu_0](dx) = \int_{\overline{\mathscr{S}^d}(\R)} \C \mu_{0, x} (\d \C) \rho(dx).
    \end{align*}
    Clearly $\int_{\overline{\mathscr{S}^d}(\R)} \C \mu_{0, x} (\d \C) \in L^2(\Omega; d \rho)$ and $\mathbb{P}[\mu_0] \in L^1(\Omega)$ if $\rho \in L^1(\Omega)$. As in \eqref{eq:PL2_ineq} we find
    \begin{align*}
        \bigg{|} \int_{\overline{\mathscr{S}^d}(\R)} \C \left( \mu_{0,x}(\d \C) - \mu_{1, x}(\d \C) \right) \bigg{|} \le W_1 \left( \mu_{0,x}, \mu_{1,x} \right) \le W_2 \left( \mu_{0,x}, \mu_{1,x} \right)
    \end{align*}
    which implies
    \begin{align*}
        \| \mathbb{P}[\mu_0] - \mathbb{P}[\mu_1] \|_{L^1(\Omega)} \le d_{Q_{\rho}} \left( \mu_0, \mu_1 \right).
    \end{align*}
    We conclude that $E : D_{Q_{\rho}} \rightarrow \R^+$ is lower semicontinuous if $\rho < < \mathscr{L}^d_{\Omega}$, where $\mathscr{L}^d_{\Omega}$ is the $d$-dimensional Lebesuge measure on $\Omega$.

    Now, take the initial datum $\mu(t=0) = \mu(x, \C; t=0)$ of the reduced Wasserstein gradient flow \eqref{eq:wasstype_mu} in the leaf $Q_{\rho}$, with disintegration $\mu_x (t=0) = \mu_x(\C; t=0)$ (with respect to $\rho$). Then, the solution $\mu(t)$ remains in $Q_{\rho}$ (as long as it exists) and its disintegration $\mu_x( \C ; t)$ with respect to $\rho$ satisfies the same flow equation \eqref{eq:wasstype_mu} $\rho$-a.e., with initial datum $\mu_x(t=0)$. The Poisson equation for the pressure $p=p(x,t)$ now reads
    \begin{align*}
        - \diver \left( \left( \rho \mathbb{P}[\mu_x] + r \mathbb{I} \right) \grad p \right) = S \quad \mathrm{in} \; \Omega.
    \end{align*}
\end{remark}

\subsection{Scalar model.} A modification of the scalar mesoscopic model, as proposed in \cite{burger2019}, is obtained by imposing the following structure for the measure $\mu$,
\begin{align*}
    \mu(x, \C) = \int_0^{\infty} \int_{\mathbb{S}^1_+} \delta \left( \C - C \theta \otimes \theta \right) \eta(x, \d \theta, \d C),
\end{align*}
where $C = C(x, \theta) \in \R^+_0$ represents the scalar conductivity of the network at the point $x$ in direction $\theta \in \mathbb{S}^1_+ := \{ v \in \R^d , \mathrm{s.t.} , |v|=1, v_1 \ge 0 \}$, and $\eta$ is a probability measure on $\overline{\Omega} \times \mathbb{S}^1_+ \times \R^+$. Consequently, taking $r: = 0$ in accordance with \cite{burger2019}, the energy functional \eqref{eq:energy} takes the form
\begin{align*}
    \tilde{E}[\eta] = \int_{\R^+} \int_{\mathbb{S}^1_+} \int_{\Omega} \left(  C | \theta \cdot \grad p[\eta] |^2 + \frac{\nu}{\gamma} C^{\gamma} \right) \d \eta,
\end{align*}
where we interpret $\eta(x, \theta, C)$ as the probability of having an edge with a conductivity $C \ge 0$ in the direction $\theta$ at the position $x \in \Omega$. Here, $p[\eta]$ solves the Poisson equation \eqref{eq:poisson1}-\eqref{eq:bc_poisson} with the permeability tensor given by
\begin{align}
\label{eq:burger_permeability}
    \tilde{\P} [\eta] (x) = \int_{ \R_+} \int_{\mathbb{S}^1_+} C \theta \otimes \theta \eta( x, \d \theta, \d C).
\end{align}
Proceeding as in Proposition \ref{wasstype_GF}, we derive a reduced Wasserstein gradient flow for $\eta$,
\begin{align}
    \label{eq:burger_mu}
    \partial_t \eta + \partial_C \left( \eta \left(  | \theta \cdot \grad p|^2 - \nu C^{\gamma-1} \right) \right) = 0.
\end{align}

\section{Stationary solutions}\label{sec:stationary}

In this section we explore stationary solutions of the PDE systems introduced in Section \ref{sec:meso}.

\subsection{Reduced Wasserstein model} \label{subsec:reduced_wass}
From the stationary version of \eqref{eq:wasstype_mu} we immediately recollect that the stationary measure of the system $\mu_{\infty}$ is concentrated on the set
\begin{align*}
     \big{ \{ } (x, \C) \in X : \grad p \otimes \grad p = \nu | \C |^{\gamma-2} \C \big{ \} }.
\end{align*}
Consequently, for $\gamma>1$,
\begin{align}
    \label{eq:stationary_wasstype}
    \mu_{\infty} (x, \C) = \rho(x) \delta \left( \C - k\left( | \grad p | \right) \grad p \otimes \grad p \right).
\end{align}
Note that $\mu_{\infty} \in Q$ if and only if $\rho(x) = 1$ a. e. in $\Omega$, but we can allow for more mathematical generality here, by choosing $\rho \in \mathscr{P}(\overline{\Omega})$, referring to Remark \ref{remark2}.
Indeed, plugging \eqref{eq:stationary_wasstype} into \eqref{eq:wasstype_mu} we obtain
\begin{align}
     \grad p \otimes \grad p = \nu k \left( | \grad p | \right)^{\gamma - 1} | \grad p |^{2(\gamma - 2)} \grad p \otimes \grad p, \notag
\end{align}
which gives
\begin{align}
    \label{eq:k_gradp}
    k \left( | \grad p | \right) = \left(  \frac{1}{\nu | \grad p |^{2(\gamma - 2)}} \right)^{\frac{1}{\gamma - 1}}.
\end{align}
Substituting \eqref{eq:stationary_wasstype} into the permeability tensor \eqref{eq:permeability_tensor} yields
\begin{align*}
    \P [\mu_{\infty}] (x) = \rho(x) k \left( | \grad p | \right) \grad p \otimes \grad p,
\end{align*}
and the Poisson equation reads
\begin{align}
    \label{eq:stationary_Poisson}
    - \diver \bigg[ \left( r + \left( \frac{1}{\nu} \right)^{\frac{1}{\gamma - 1}} \rho(x) | \grad p |^{\frac{2}{\gamma - 1}} \right) \grad p \bigg] = S.
\end{align}
subject to the zero-flux boundary condition on $\partial \Omega$. Note that $\rho$ needs to be absolutely continuous with respect to the $d$-dimensional Lebesgue measure on the set $\{ x \in \Omega : \grad p(x) \neq 0 \}$ (see the comment below \eqref{eq:domain_E}) and obviously the Poisson equation \eqref{eq:stationary_Poisson} does not 'see' any Lebesgue-singular part of $\rho$.

\begin{theorem}\label{stationary_wasstype_gamma>1}
    For any $S \in L^2(\Omega)$, $\rho \in L^{\infty} (\Omega)$ and $\gamma > 1$ there exists a unique solution $p \in \overline{H}^1(\Omega)$, $ \grad p \in L^{\frac{2 \gamma}{\gamma - 1}}(\Omega; \d \rho)$ of \eqref{eq:stationary_Poisson} subject to the zero-flux boundary condition on $\partial \Omega$.
\end{theorem}

The space $\overline{H}^1(\Omega)$ is the space of all $H^1(\Omega)$ functions with zero integral mean. Notice that for $\gamma > 1$ we have $\frac{2 \gamma}{\gamma - 1} > 2$.
    
\begin{proof}
    Let us define the functional $\mathcal{J} : \overline{H}^1(\Omega) \rightarrow \R \cup \{ \infty \}$ as
    \begin{eqnarray}
    \label{eq:functional1}
        {\mathcal{J}[p] :=} \left\{ 
        \begin{array}{ll} \displaystyle
        \bigintsss_{\Omega} \bigg[ \frac{r}{2} | \grad p |^2 + \frac{\gamma - 1}{2 \gamma} \left( \frac{1}{\nu} \right)^{\frac{1}{\gamma-1}} \rho(x) | \grad p |^{\frac{2 \gamma}{\gamma - 1}} - p S \bigg] \d x & \mathrm{if} \, \grad p \in L^{\frac{2 \gamma}{\gamma - 1}}(\Omega; \d \rho), \\ \displaystyle
        + \infty & \mathrm{if} \, \grad p \not\in L^{\frac{2 \gamma}{\gamma - 1}}(\Omega; \d \rho),
        \end{array}
        \right.
    \end{eqnarray}
    with associated Lagrangian
    \begin{align*}
        \mathcal{L}\left( q, p, x \right) := \frac{r}{2} | q |^2 + \frac{\gamma - 1}{2 \gamma} \left( \frac{1}{\nu} \right)^{\frac{1}{\gamma-1}} \rho(x) | q |^{\frac{2 \gamma}{\gamma - 1}} - p S.
    \end{align*}
    For simplicity of notation let us consider ${\nu} = 1$. Then uniform convexity is guaranteed since
    \begin{align*}
        \mathcal{L}_{q_{i} q_{j}} &= \frac{\partial}{\partial q_{j}} \left( r
        q_{i} + \rho(x) q_{i} | q |^{\frac{2}{\gamma-1}} \right) \\
        &= \left( r + \rho(x) | q |^{\frac{2}{\gamma - 1}} + \frac{2}{\gamma - 1} \rho(x) q_{i} q_{j} | q |^{\frac{2(2 - \gamma)}{\gamma - 1}} \right) \delta_{i j},
    \end{align*}
    and for every $ \xi \in \R^d$,
    \begin{align*}
        \sum_{i, j = 1}^d \mathcal{L}_{q_{i} q_{j}} \xi_i \xi_j = \left( r + \rho(x) \left( 1 + \frac{2}{\gamma-1} \right) | q |^{\frac{2}{\gamma - 1}} \right) | \xi |^2 \ge r | \xi |^2,
    \end{align*}
    for $\gamma > 1$. Coercivity can be deduced from an application of Hölder inequality and Poincaré-Wirtinger inequality on $\overline{H}^1 (\Omega)$. As a result, the existence of a unique minimizer for \eqref{eq:functional1} is guaranteed by the standard theory of calculus of variations \cite{evans}. Finally, we notice that \eqref{eq:stationary_Poisson} represents the Euler-Lagrange equation associated with the minimization problem \eqref{eq:functional1}.
    \end{proof}

We underline that Theorem \ref{stationary_wasstype_gamma>1} holds if $r=0$ and $\rho$ uniformly positive, while in the stated version $r>0$ and $\rho$ may attain zero. 
Note also that \eqref{eq:stationary_Poisson} is the $p$-Laplace equation with $p = \frac{2\gamma}{\gamma-1}$,
see, e.g., the works \cite{ladyzhenskaya, nazarov, urdaletova} of N. N. Uraltseva.

For $\gamma=1$ the stationary version of \eqref{eq:wasstype_mu} imposes that $\mu_{\infty}$ be concentrated on the set
\begin{align*}
     \left\{ (x, \C) \in X : \grad p \otimes \grad p = \nu \frac{\C}{|\C|} \right\},
\end{align*}
where for $\C=0$ we interpret the expression ${\C}/{|\C|}$ as the subdifferential of $|\C|$, i.e., the closed unit ball $\{|\C|\leq 1\}$.
Consequently, we have $|\grad p|^2 \leq \nu$, and $|\grad p|^2=\nu$ whenever $\C\neq 0$.
Therefore, there exists a positive measurable function $\lambda=\lambda(x)$ such that
\begin{align}
    \label{eq:stationary_wasstype1}
    \mu_{\infty} (x, \C) = \rho(x) \delta \left( \C - \lambda \chi_{|\grad p|=\sqrt\nu} \grad p \otimes \grad p \right),
\end{align}
where $\rho \in \mathscr{P}(\overline{\Omega})$ and 
$\chi_{|\grad p|=\sqrt\nu}$ denotes the characteristic function of the set $\{|\grad p|=\sqrt\nu\}$.

Substituting \eqref{eq:stationary_wasstype1} into the permeability tensor \eqref{eq:permeability_tensor} yields
\begin{align*}
    \P [\mu_{\infty}] (x) = \rho(x) \lambda(x) \chi_{|\grad p|=\sqrt\nu} \grad p \otimes \grad p,
\end{align*}
and we arrive at the highly nonlinear Poisson equation
\begin{align}
    \label{eq:stationary_Poisson1}
    - \diver \bigg[ \left( r + \nu \rho(x) \lambda(x) \chi_{|\grad p|=\sqrt\nu} \right) \grad p \bigg] = S
\end{align}
subject to the zero-flux boundary condition on $\partial \Omega$.
An analogous problem was studied in \cite[Section 4.2]{HMPS16}, where it was shown 
to be equivalent to the free boundary problem
\begin{align}
   -\grad\cdot\left[(r+\rho(x)a(x)^2)\grad p\right] &= S,\qquad p \in H^1_0(\Omega), \label{NLP1} \\
    |\grad p(x)|^2 &\leq \nu,\qquad \mbox{a.e. on }\Omega, \label{NLP2} \\
   \rho(x) a(x)^2 \left[|\grad p(x)|^2-\nu \right] &= 0,\qquad \mbox{a.e. on }\Omega, \label{NLP3}
\end{align}
for some measurable function $a^2=a(x)^2$ on $\Omega$ which is the Lagrange multiplier
for the condition \eqref{NLP2}.
The function $\lambda=\lambda(x)$ can be chosen as $\lambda(x):={a(x)^2}/{\nu}$.
Solutions of the system \eqref{NLP1}--\eqref{NLP3} were constructed in \cite[Lemma 5]{HMPS16} as the unique minimizers of the convex and coercive functional
\begin{equation} \label{eq:twisting}
    \mathcal{I}[p] := \int_\Omega \left( r \frac{|\grad p|^2}{2} - p S \right) \d x
\end{equation}
on the set $\{p\in H_0^1(\Omega), |\grad p|^2\leq \nu \mbox{ a.e. on } \Omega\}$.
The gradient constrained variational problem \eqref{eq:twisting} was studied in \cite{caff} as a model for twisting of an elastic–plastic cylindrical bar. There it was shown that the unique solution has $C^{1,1}$-regularity in $\Omega$, see
also \cite{Safdari1, Safdari2}.
In this context let us point out the significant works \cite{apushkinskaya, shago, uraltseva} of N. N. Uraltseva on regularity and other fine properties of solutions of free boundary problems.

\subsection{Wasserstein model.} Studying the stationary case of \eqref{eq:wass_mu} for $\gamma>1$ we notice that $\mu_{\infty}$ is concentrated on the set
\begin{align*}
    &\big{\{} (x, \C) \in X : \grad p \otimes \grad p = \nu | \C |^{\gamma-2} \C, \; \grad\left(\grad p \cdot \C \grad p \right) = 0 \big{\}} \\
    =& \big{\{} (x, \C) \in X : \C = k(| \grad p |) \grad p \otimes \grad p, \; \grad p(x) = 0 \; \mathrm{or} \; \grad | \grad p(x) | = 0 \big{\}},
\end{align*}
where $k(| \grad p |)$ is defined as in \eqref{eq:k_gradp}.
Consider the following subset of $\overline{\Omega}$
\begin{align*}
    \Lambda := \big{\{} x \in \overline{\Omega} : \grad | \grad p(x) | = 0, \; \grad p(x) \neq 0 \big{\}}
\end{align*}
we conclude that $\mu_{\infty}$ is concentrated on
\begin{align*}
    &\big{\{} (x, \C) \in X : \C = k(|\grad p|) \grad p \otimes \grad p, \, x \in \Lambda \big{\}} \cup \{ (x, \C) \in X : \grad p(x) = 0 \; \mathrm{and} \; \C = 0 \}.
\end{align*}
As in Section \ref{subsec:reduced_wass}, we conclude
\begin{align}
    \label{eq:stationary_wass}
    \mu_{\infty} (x, \C) = \rho(x) \delta \left( \C - k(| \grad p |) \grad p \otimes \grad p \right),
\end{align}
where a possible Lebesgue singular part of $\rho$ concentrates on $\{ x \in \overline{\Omega} : \grad p(x) = 0 \}$. Again, the latter does not show up in the permeability tensor. 
If $p \in H^2(\Omega)$ then $| \grad p | \in H^1(\Omega)$, see \cite[Theorem 6.17]{lieb} and $\grad | \grad p | = 0$ holds Lebesgue almost everywhere on the set $\{ x \in \Omega : \grad p(x) = 0 \}$, see \cite[Theorem 6.19]{lieb}. This implies that the Lebesgue absolutely continuous part of $\rho$ concentrates in the set $\Lambda$.

Now, assume $\rho \in C^1(\overline{\Omega})$ and $\mathrm{supp} \rho = \bigcup_{i=1}^N \Gamma_i$, where $N$ may be infinite and $\Gamma_i$ are the pairwise disjoint connected components. Consequently, there exist constants $\omega_i$ such that $| \grad p | = \omega_i$ in $\Gamma_i$. Then, $p$ has to satisfy
\begin{align}
    \label{eq:stationary_poisson_wass}
    - \diver \left( \left( r + \rho(x) \sum_{i=1}^N \mathbbm{1}_{\Gamma_i} \left( \frac{\omega_i^2}{\nu} \right)^{\frac{1}{\gamma-1}} \right) \grad p \right) = S
\end{align}
subject to zero-flux boundary condition on $\partial \Omega$. This leads to the following Theorem.

\begin{theorem} \label{thm:stationary_wass}
    Let $S \in L^2(\Omega)$ with zero integral mean and let $\gamma > 1$. Then, \eqref{eq:stationary_wass} constitutes an equilibrium solution of \eqref{eq:wass_mu} if and only if there are constants $\omega_1^2, \cdots, \omega_N^2 $ such that the solution $p$ of \eqref{eq:stationary_poisson_wass} with zero-flux boundary conditions satisfies
    \begin{align*}
        | \grad p (x) | = \omega_i \quad \mathrm{on} \; \Gamma_i.
    \end{align*}
\end{theorem}

Note that the only possible equilibrium solutions are of the form given in Theorem \ref{thm:stationary_wass}. 
Although these equilibrium solutions can be called 'non-generic', it is actually an easy exercise to construct examples.
Indeed, take a smooth function $p=p(x)$ on $\overline{\Omega}$ such that $ | \grad p | = \mathrm{const.}$ on a subdomain $\Sigma \subset \subset \Omega$ and $\grad p \cdot n = 0$ on $\partial \Omega$.
Then, choose a smooth probability density $\rho \in \mathscr{P}(\overline{\Omega})$ which concentrates in $\Sigma$ and compute $S = S(x)$ from the Poisson equation \eqref{eq:stationary_Poisson}.

\subsection{Monokinetic model}
The stationary solutions for the monokinetic model \eqref{eq:monokinetic_wass_rho}-\eqref{eq:monokinetic_L2_C} are a mix of the Wasserstein and the reduced Wasserstein model. Indeed, taking the stationary version of \eqref{eq:monokinetic_L2_C} we conclude
\begin{align}
    \hat{\C}(x) = k \left( |\grad p| \right) \grad p \otimes \grad p, \; \; \mathrm{or} \; \; \rho(x) \equiv 0 \notag
\end{align}
with $k \left( |\grad p| \right) = \left( \frac{1}{\nu | \grad p |^{2(\gamma - 2)}}\right)^{\frac{1}{\gamma-1}}$, in analogy with Section \ref{subsec:reduced_wass}. Then, the Poisson equation becomes \eqref{eq:stationary_Poisson} (note that only the product of $\rho$ and $\C$ appear in the Poisson equation). 

On the other hand, considering the stationary version of \eqref{eq:monokinetic_wass_rho} coupled to the above $\hat{\C}$ leads to
\begin{align*}
    \grad \bigg[ k \left( | \grad p | \right) | \grad p |^4 - \frac{\nu}{\gamma} k\left( | \grad p | \right)^{\gamma} | \grad p |^{2\gamma} \bigg] = \left( \frac{1}{\nu} \right)^{\frac{1}{\gamma-1}} \left( \frac{\gamma - 1}{\gamma} \right) \grad \bigg[ | \grad p |^{\frac{2\gamma}{\gamma-1}}\bigg] = 0.
\end{align*}
If $\gamma \neq 1$, $|\grad p|$ will be constant. Consequently, the latter scenario is identical to the stationary Wasserstein model \ref{subsec:GFWass}, and Theorem \ref{thm:stationary_wass} applies.

\subsection{Scalar model} 
From the stationary version of \eqref{eq:burger_mu} we immediately notice that
\begin{align*}
    \left(  | \theta \cdot \grad p |^2 - \nu C^{\gamma-1} \right) \eta_{\infty} = 0,
\end{align*}
which implies that $\eta_{\infty}$ is concentrated on
\begin{align*}
    \big{ \{ } (x, \theta, C) \in \overline{\Omega} \times \mathbb{S}^1_+ \times \R^{+} : | \theta \cdot \grad p |^2 = \nu C^{\gamma-1} \big{ \} }.  
\end{align*}
Consequently,
\begin{align}
    \label{eq:burger_stationary}
    \eta_{\infty} = \rho(x, \theta) \delta \left( C - \left( \frac{1}{\nu} \right)^{\frac{1}{\gamma-1}} | \grad p \cdot \theta|^{\frac{2}{\gamma - 1}}  \right),
\end{align}
where $\rho \in \mathscr{P}^{+}\left(\overline{\Omega} \times \mathbb{S}^1_+ \right)$. Then, the permeability tensor reads
\begin{align*}
    \tilde{\mathbb{P}}[\eta_{\infty}] (x) = \left( \frac{1}{\nu} \right)^{\frac{1}{\gamma - 1}} \int_{\mathbb{S}^1_+} | \grad p \cdot \theta |^{\frac{2}{\gamma - 1}} \theta \otimes \theta \rho(x, \d \theta),
\end{align*}
and the Poisson equation is given by
\begin{align}
    \label{eq:stationary_Poisson_burger}
    - \left( \frac{1}{\nu} \right)^{\frac{1}{\gamma - 1}} \diver \bigg[ \left( \int_{\mathbb{S}^1_+} | \grad p \cdot \theta |^{\frac{2}{\gamma - 1}} \theta \otimes \theta \rho (x , \d \theta) \right) \grad p \bigg] = S.
\end{align}
We reiterate that here we set $r=0$ in accordance with \cite{burger2019}.

\begin{theorem}
    For any $S \in L^2(\Omega)$, $\rho \in L^{\infty} (\Omega \times \mathbb{S}^1_+)$ uniformly positive and $\gamma > 1$ there exists a unique solution $p \in \overline{H}^1(\Omega) \cap W^{1, \frac{2 \gamma}{\gamma - 1}}(\Omega)$ of \eqref{eq:stationary_Poisson_burger} subject to the zero-flux boundary condition on $\partial \Omega$.
\end{theorem}

\begin{proof}
    The proof mirrors the one of Theorem \ref{stationary_wasstype_gamma>1}, with \eqref{eq:stationary_Poisson_burger} serving as the Euler-Lagrange equation for the functional $\mathcal{J} : \overline{H}^1(\Omega) \rightarrow \R \cup { \infty }$, defined as follows
    \begin{eqnarray}
    \label{eq:functional_J}
        {\mathcal{J}[p] :=} \left\{ 
        \begin{array}{ll} \displaystyle
        \bigintssss_{\Omega} \bigg[ \int_{\mathbb{S}^1_+} \left( \frac{1}{\nu} \right)^{\frac{1}{\gamma - 1}} \frac{\gamma - 1}{2 \gamma}  \rho(x, \theta) | \theta \cdot \grad p |^{\frac{2 \gamma}{\gamma - 1}} \d \theta - p S \bigg] \d x & \mathrm{if} \, \grad p \in L^{\frac{2 \gamma}{\gamma - 1}}(\Omega), \\ \displaystyle
        + \infty & \mathrm{if} \, \grad p \not\in L^{\frac{2 \gamma}{\gamma - 1}}(\Omega).
        \end{array}
        \right.
    \end{eqnarray}
\end{proof}

\subsection{Semi-discrete model}
Now we take a relaxation of the functional \eqref{eq:functional_J}, that is
\begin{eqnarray}
\label{eq:functional_J2}
    {\mathcal{\tilde{J}}[p] :=} \left\{ 
    \begin{array}{ll} \displaystyle
    \left( \frac{1}{\nu} \right)^{\frac{1}{\gamma - 1}} \frac{\gamma - 1}{2 \gamma} \bigintssss_{\overline{\Omega}} \int_{\mathbb{S}^1_+} | \theta \cdot \grad p |^{\frac{2 \gamma}{\gamma - 1}} \d \rho - \int_{\Omega} p S \d x & \mathrm{if} \, \theta \cdot \grad p \in L^{\frac{2 \gamma}{\gamma - 1}}(\overline{\Omega} \times \mathbb{S}^1_+, d \rho), \\ \displaystyle
    + \infty & \mathrm{if} \, \theta \cdot \grad p \not\in L^{\frac{2 \gamma}{\gamma - 1}}(\overline{\Omega} \times \mathbb{S}^1_+, d \rho).
    \end{array}
    \right.
\end{eqnarray}
This functional allows us to obtain new (semi-)discrete transportation models on graphs in $\overline{\Omega}$. Therefore we think of the graph $\mathcal{G}$, embedded into $\overline{\Omega}$, as consisting of the finite vertex set $\mathcal{V} = \{ x_1, \dots, x_N \} \subseteq \overline{\Omega}$, connected by smooth open edges
\begin{align*}
    \mathcal{E} := \big{\{} e_{ij} : (0, L_{ij}) \rightarrow \overline{\Omega} : (i,j) \subseteq \{ 1, \dots, N \}^2, \; i< j, \; e_{ij}(0^+) = x_i, \; e_{ij}(L_{ij}^-) = x_j \big{\}}.
\end{align*}
For the sake of simplicity we assume that $e_{ij} = e_{ij}(s)$ is parametrized by the arclength parameter $s \in (0, L_{ij})$, such that $L_{ij}$ is the length of the edge connecting $x_i$ to $x_j$. The (open) edges are not allowed to intersect and always connect different vertices. Indeed, note that if $e_{ij} \in \mathcal{E}$ then $e_{ij} \not\in \mathcal{E}$ such that each pair of vertices $(x_i, x_j)$ can only be connected by at most one edge. We denote by $t_{ij} = t_{ij}(x)$ the unit tangent vector of $e_{ij}$ at the point $x = e_{ij}(s)$ such that $t_{ij}(x) = \frac{d}{ds} e_{ij}(s)$.

Now let $\mathscr{H}^1$ be the one-dimensional Hausdorff measure on $\R^d$ and $\mathscr{H}_{ij}^1$ its restriction to the open edge $e_{ij}$, i.e.,
\begin{align*}
    \mathscr{H}_{ij}^1 (A) = \mathscr{H}^1 \left( A \cap \{ e_{ij} : 0 < s < L_{ij} \}\right)
\end{align*}
for every Borel set $A \subseteq \R^d$. Let $\beta : \mathcal{G} \rightarrow \R^+$ be a bounded positive function defined on the (closed) graph $\mathcal{G}$ with $\sum_{e_{ij} \in \mathcal{E}} \int_{e_{ij}} \beta \d s = 1$ and define the probability measure $\rho \in \mathscr{P} \left( \overline{\Omega} \times \mathbb{S}_+^1 \right)$ by
\begin{align*}
    \rho(x, \theta) = \sum_{e_{ij} \in \mathcal{E}} \beta(x) \mathscr{H}_{ij}^1(x) \otimes \delta \left( \theta - t_{ij}(x) \right).
\end{align*}
For the following we shall assume that $S$ concentrates on the graph $\mathcal{G}$:
\begin{align*}
    S(x) = \sum_{i=1}^N S_i \delta(x - x_i) + \sum_{e_{ij} \in \mathcal{E}} S_{ij}(x) \mathscr{H}_{ij}^1(x),
\end{align*}
where $S_i \in \R$ and $S_{ij} : e_{ij} \rightarrow \R$ for all $i, j = 1, \dots, N$ such that the global conservation of mass \eqref{eq:cont_globalconservation} applies.

Then the functional \eqref{eq:functional_J2} becomes
\begin{eqnarray*}
    {\mathcal{\tilde{J}}[p] :=} \left\{ 
    \begin{array}{ll} \displaystyle
    \left( \frac{1}{\nu} \right)^{\frac{1}{\gamma - 1}} \frac{\gamma - 1}{2 \gamma} \sum_{e_{ij}\in \mathcal{E}} \int_{e_{ij}} \beta | \partial_s p |^{\frac{2 \gamma}{\gamma - 1}} \d s - \left( \sum_{i=1}^N S_i p(x_i) + \sum_{e_{ij} \in \mathcal{E}} \int_{e_{ij}} S_{ij}p \d s \right) & \mathrm{if} \, \partial_s p \in L^{\frac{2 \gamma}{\gamma - 1}}(e_{ij}), \\ \displaystyle
    + \infty & \mathrm{if} \, \partial_s p \not\in L^{\frac{2 \gamma}{\gamma - 1}}(e_{ij}).
    \end{array}
    \right.
\end{eqnarray*}
Note that the relaxation of the functional $J$ manifests itself in the fact that now the permeability tensor $\mathbb{P}$ is a singular measure on $\overline{\Omega}$ concentrating on the graph $\mathcal{G}$ such that the Poisson equation \eqref{eq:poisson1}-\eqref{eq:bc_poisson} for p is a-priorily  meaningless on $\overline{\Omega}$. 
However, looking for critical points of the relaxed functional bypasses the Poisson equation on $\overline{\Omega}$ and it is an easy exercise to show that the Euler-Lagrange equation of $\mathcal{\tilde{J}}[p]$ is given by
\begin{align}
    \label{eq:semidiscrete}
    - \left( \frac{1}{\nu} \right)^{\frac{1}{\gamma-1}} \frac{\partial}{\partial s}\left( \beta \partial_s p | \partial_s p |^{\frac{2}{\gamma-1}} \right) = S_{ij} (s) \quad \mathrm{on} \; e_{ij}
\end{align}
subject to the transmission conditions \cite{lions}
\begin{align}
    \label{eq:transmission}
    - \left(\frac{1}{\nu}\right)^{\frac{1}{\gamma-1}} \beta(x_i) \sum_{j \in \mathcal{N}(i)} | \partial_s p (x_i; e_{ij})|^{\frac{2}{\gamma-1}} \partial_s p (x_i; e_{ij}) = S_i \quad \forall i=1,\dots,N.
\end{align}
Here, $\partial_s p (x_i; e_{ij})$ denotes the value of the $s$-derivative of $p$ computed at the node $x_i$ along the edge $e_{ij}$. Since $p$ is supposed to be continuous on $\mathcal{G}$ but not necessarily smooth, jumps in the derivative can occur and the above introduced notation is essential to distinguish between pressure derivatives originating from distinct edges.

The semi-discrete model \eqref{eq:semidiscrete}-\eqref{eq:transmission} significantly enriches the discrete models presented in \cite{burger2019, hucai2013}.

Now, we consider the special case where $S$ is concentrated only on the graph nodes, i.e., $S_{ij} \equiv 0 \; \forall e_{ij} \in \mathcal{E}$. We set $\beta(x) = \frac{1}{L_{ij}}$ on $e_{ij}$ and \eqref{eq:semidiscrete} reads
\begin{align*}
    \frac{\partial}{\partial s} \left( \partial_s p | \partial_s p |^{\frac{2}{\gamma-1}} \right) = 0 \quad \mathrm{on} \; e_{ij},
\end{align*}
implying $\partial_s p = K_{ij}$ on $e_{ij}$ where $K_{ij} \in \R$ is a constant. Thus, $p$ is a linear function on $e_{ij}$, i.e., $p(s) = K_{ij} s + c$ on $e_{ij}$ with $c \in \R$. Note that $P_i:=p(x_i)=p(0)=c$ and $P_j:=p(x_j)=p(L_{ij})= \alpha L_{ij} + c$ so that
\begin{align*}
    \partial_s p = \frac{P_j - P_i}{L_{ij}} \quad \mathrm{on} \; e_{ij}.
\end{align*}
Substituting the latter expression into the transmission condition \eqref{eq:transmission} we obtain
\begin{align}
    \label{eq:consistency_disc_semidisc}
    - \left( \frac{1}{\nu} \right)^{\frac{1}{\gamma-1}} \frac{1}{L_{ij}} \sum_{j \in \mathcal{N}(i)} \bigg| \frac{P_j - P_i}{L_{ij}} \bigg|^{\frac{2}{\gamma-1}} \frac{P_j - P_i}{L_{ij}} = S_i \quad \mathrm{for} \; i=1,\dots,N.
\end{align}
This equation on the nodes of $\mathcal{G}$ is coherent with the discrete model outlined in Section \ref{sec:discrete} with $r=0$. Indeed, consider the stationary version of \eqref{eq:discrete_GF}
\begin{align*}
    \frac{Q_{ij}^2}{C_{ij}^2} = \nu C_{ij}^{\gamma-1}
\end{align*}
and by plugging in the flow $Q_{ij}$ \eqref{eq:discrete_flow} we get
\begin{align*}
    C_{ij} = \left( \frac{1}{\nu} \right)^{\frac{1}{\gamma-1}} \bigg| \frac{P_j - P_i}{L_{ij}} \bigg|^{\frac{2}{\gamma-1}}.
\end{align*}
Finally, by substituting the latter into the Kirchhoff law \eqref{eq:kirchhoff_Law} we find \eqref{eq:consistency_disc_semidisc}.

\section{Fisher-Rao gradient flow}\label{sec:fisher_rao}
An important metric on the probability space $\mathscr{P}(X)$ is the Fisher-Rao metric - sometimes referred to as the Hellinger-Kakutani distance (see \cite{liero2016, liero2018}) - which arises naturally in information geometry \cite{rao}. Our study is restricted to the space of smooth, positive, and absolutely continuous probability measures $\mu$, denoted by $\mathrm{Prob}(X)$. Note that the Fisher-Rao metric is -up to a multiplicative factor and for manifolds of dimension $>1$- the unique smooth metric invariant under the action of the diffeomorphism group on the manifold, see \cite{bauer}.

We define the gradient structure triple $\left( \mathrm{Prob}(X), E, \mathbb{G}^{\mathrm{FR}} \right)$, where $E$ is defined by \eqref{eq:energy} coupled to the Poisson equation \eqref{eq:poisson1} with no-flux boundary condition and its domain is $\big{\{} \mu \in \mathrm{Prob}(X) : \int_{X} \left( | \C| + |\C|^{\gamma} \right) \d \mu < \infty \big{\}}$. Moreover, $\mathbb{G}^{\mathrm{FR}}(\mu) : T_{\mu} \mathrm{Prob}(X) \rightarrow T_{\mu}^* \mathrm{Prob}(X)$, $\mathbb{G}^{\mathrm{FR}}(\mu) v := \frac{v}{\mu}$ for every $v \in T_{\mu} \mathrm{Prob}(X)$.
Naturally, $\mathbb{G}^{\mathrm{FR}}(\mu)$ induces the metric tensor
\begin{align*}
    g_{\mu}^{\mathrm{FR}}(v, \tilde{v}) := \leftindex_{T_{\mu}^* \mathrm{Prob}(X)} {\braket{ \mathbb{G}^{\mathrm{FR}}(\mu) v, \tilde{v}}}_{T_{\mu} \mathrm{Prob}(X)} = \int_{X} \frac{v \tilde{v}}{\mu^2} \d \mu \quad \forall \, v, \tilde{v} \in T_{\mu} \mathrm{Prob}(X).
\end{align*}
Note that the metric tensor is precisely the Fisher information of information science \cite{rao}. Moreover, such metric tensor induces the Fisher-Rao (or spherical Hellinger-Kakutani) distance (see \cite{friedrich, laschos-mielke, liero-mielke-sav})
\begin{align}
    \label{eq:fisher-rao_distance}
    d_{\mathrm{FR}}(\mu_0, \mu_1) = 2 \arccos{\int_{X} \sqrt{\mu_0} \sqrt{\mu_1} \d x \d \C}, \; \forall \mu_0, \mu_1 \in \mathscr{P}(X).
\end{align}
Next, we introduce the operator $\mathbb{K}^{\mathrm{FR}}(\mu) : T_{\mu}^* \mathrm{Prob}(X) \rightarrow T_{\mu} \mathrm{Prob}(X)$ as the inverse operator of $\mathbb{G}^{\mathrm{FR}}(\mu)$, i.e., $\mathbb{K}^{\mathrm{FR}}(\mu) := \mathbb{G}^{\mathrm{FR}}(\mu)$. This operator is defined via \cite{wang2020}
\begin{align*}
    \mathbb{K}^{\mathrm{FR}}(\mu) \Phi = \mu \left( \Phi - \int_X \Phi \d \mu \right) \quad \forall \Phi \in T_{\mu}^{*} \mathrm{Prob}(X). 
\end{align*}
Consequently, the associated gradient flow equation reads
\begin{align*}
    \partial_t \mu = - \mathbb{K}^{\mathrm{FR}}(\mu) \frac{\delta E[\mu]}{\delta \mu} = - \mu \left( \frac{\delta E[\mu]}{\delta \mu} - \int_{X} \frac{\delta E[\mu]}{\delta \mu} \d \mu \right).
\end{align*}
In our context, plugging the first variation \eqref{eq:first_variation} in the latter equation, we obtain
\begin{align}
    \label{eq:fr_gradientflow}
    \partial_t \mu = \mu \left( \grad p[\mu] \cdot \C \grad p[\mu] - \frac{\nu}{\gamma} | \C |^{\gamma} - \int_{X} \left( \grad p[\mu] \cdot \C \grad p[\mu] - \frac{\nu}{\gamma} | \C |^{\gamma}\right) \d \mu \right),
\end{align}
coupled to the Poisson equation \eqref{eq:poisson1}-\eqref{eq:bc_poisson}.

The energy dissipation relation reads
\begin{align*}
    \frac{\d}{\d t} E[\mu(t)] = - \bigintsss_{X} \left| \grad p \cdot \C \grad p - \frac{\nu}{\gamma} |\C|^{\gamma} \right|^2 \d\mu
    + \left( \bigintsss_{X} \left(\grad p \cdot \C \grad p - \frac{\nu}{\gamma} | \C |^{\gamma}\right) \d\mu \right)^2 \leq 0.
\end{align*}

We now analyze the lower semicontinuity of $E[\mu]$ with respect to the Fisher-Rao metric \eqref{eq:fisher-rao_distance}.
\begin{theorem} \label{thm:fisher-rao_lsc}
    Let $\mu_n, \mu \in \mathscr{P}_{\mathrm{a.c.}}(X) := \{ \mu \in \mathscr{P}(X) : \mu \; \mathrm{abs. \; cont.} \; \mathrm{w. \; r. \; t. \; Lebesgue \; measure} \; d x d \C \; \mathrm{on} \; X \}$, $d_{\mathrm{FR}}(\mu_n, \mu) \rightarrow 0$ as $n \rightarrow \infty$ and assume that there exists $\beta > 1$ such that $| \C |^{\beta} \mu_n \in L^1(\Omega)$ uniformly for $n \in \N$. Then, $E[\mu] \le \liminf_{n \rightarrow \infty} E[\mu_n]$.
    \begin{proof}
        We immediately see that
        \begin{align*}
            \int_{X} \left( \sqrt{\mu_n} - \sqrt{\mu} \right)^2 \d x \d \C = 2 - 2 \int_{X} \sqrt{\mu_n} \sqrt{\mu} \d x \d \C \rightarrow 0 \; \mathrm{as} \; n \rightarrow \infty
        \end{align*}
        because $d_{\mathrm{FR}}(\mu_n, \mu) \rightarrow 0$ means $\int_{X} \sqrt{\mu_n} \sqrt{\mu} \d x \d \C \rightarrow 1 $ as $n \rightarrow \infty$. Thus $\sqrt{\mu_n} \rightarrow \sqrt{\mu}$ in $L^2(X)$ and $\mu_n \rightarrow \mu$ in $L^1(X)$ as well as narrowly as probability measures. After extraction of a subsequence we have $p_n := p[\mu_n] \weakconv p$ in $H^1(\Omega)$ and \eqref{eq:ubdd} gives
        \begin{align*}
            \int_{X} \big| \sqrt{\C} \grad p_n \sqrt{\mu_n} \big|^2 \d x \d \C \le K,
        \end{align*}
        uniformly in $n$, and we have $\sqrt{\C} \grad p_n \sqrt{\mu_n} \weakconv \sqrt{\C} \grad p \sqrt{\mu}$ in $L^2(X)$. This, together with the strong convergence of $\sqrt{\mu_n} \rightarrow \sqrt{\mu}$ in $L^2(X)$ implies that $\sqrt{\C} \grad p_n \mu_n \weakconv \sqrt{\C} \grad p \mu$ in $L^1(X)$. The result follows as in the proof of Theorem \ref{thm:lsc_D(E)}, using the uniform integrability of $|\C|$ with respect to the sequence $\mu_n$.
    \end{proof}
\end{theorem}
Theorem \ref{thm:fisher-rao_lsc} implies that $E$ is lower semicontinuous on $\mathscr{P}_{\mathrm{a.c.}}(X)$ with respect to the Fisher-Rao metric if $\gamma > 1$.

\subsection{Stationary Fisher-Rao}\label{subsec:stationary_fr}
From the stationary version of \eqref{eq:fr_gradientflow} we immediately recollect that every stationary measure $\mu_{\infty}$ of the system is concentrated in the set
\begin{align*}
    S := \bigg{\{} (x, \C) \in X : \frac{\delta E[\mu_{\infty}]}{\delta \mu} = \int_{X} \frac{\delta E[\mu_{\infty}]}{\delta \mu} \d \mu_{\infty} \bigg{\}}.
\end{align*}
Let $\mu_{\infty} \in \mathscr{P}(X)$, $p_{\infty} = p[\mu_{\infty}]$ and $K$ be a real constant in the range of $\frac{\delta E[\mu_{\infty}]}{\delta \mu}$, i.e., $K \in \mathrm{range}\left( \frac{\delta E[\mu_{\infty}]}{\delta \mu} \right)$, and let $\mu_{\infty}$ concentrate in
\begin{align*}
    S_K := \bigg{\{} (x, \C) \in X : \frac{\delta E[\mu_{\infty}]}{\delta \mu} = K \bigg{\}}.
\end{align*}
Then,
\begin{align*}
    \int_{X} \frac{\delta E[\mu_{\infty}]}{\delta \mu} \d \mu_{\infty} = K,
\end{align*}
and we conclude that $\mu_{\infty}$ is a stationary state. Since $\int_{X} \frac{\delta E[\mu]}{\delta \mu} \d \mu$ is in the range of $\frac{\delta E[\mu]}{\delta \mu}$ we have found all stationary states.
This leads us naturally to analyze the set, for $\mu \in \mathscr{P}(X)$ given:
\begin{align*}
    S_K = \bigg{\{} (x, \C) \in X : - \grad p[\mu] \cdot \C \grad p[\mu] + \frac{\nu}{\gamma} | \C |^{\gamma} = K \bigg{\}}.
\end{align*}
Setting $y := \grad p[\mu]$, we are lead to consider the equation
\begin{align}
    \label{eq:algebraic_equation_fr}
    \frac{\nu}{\gamma} | \C |^{\gamma} - y \cdot \C y - K = 0, \quad \mathrm{for} \; y \in \R^d, \; \C \in \overline{\mathscr{S}^d}(\R)
\end{align}
for $ K \in \mathrm{range}\left( \frac{\delta E[\mu]}{\delta \mu} \right)$.
We introduce the set
\begin{align*}
    M_K := \big{\{} (y, \C) \subset \R^d \times \overline{\mathscr{S}^d}(\R) \; \mathrm{solutions} \; \mathrm{of} \; \eqref{eq:algebraic_equation_fr} \big{\}},
\end{align*}
and prove the following propositions:

\begin{proposition}
    Let $(\mu_{\infty}, p_{\infty})$ be a stationary solution. Then, there exists $ K \in \mathrm{range}\left( \frac{\delta E[\mu_{\infty}]}{\delta \mu} \right)$ such that the probability measure $\omega := \left( \grad p_{\infty}, \mathbb{I} \right)_{\#} \mu_{\infty}$ concentrates in $M_K$.
    \begin{proof}
        Choose any two smooth test functions $\varphi=\varphi(y)$, $\psi = \psi (\C)$ and compute
        \begin{align*}
            \int_{\R^d} \int_{\overline{\mathscr{S}^d}(\R)} \varphi(y) \psi(\C) \left( \frac{\nu}{\gamma} | \C |^{\gamma} - y \cdot \C y - K \right) \d \omega &= \int_{X} \varphi(\grad p_{\infty}) \psi(\C) \left( \frac{\nu}{\gamma} | \C |^{\gamma} - \grad p_{\infty} \cdot \C \grad p_{\infty} - K \right) \d \mu_{\infty} \notag = 0
        \end{align*}
        where $ K = \int_{X} \frac{\delta E[\mu_{\infty}]}{\delta \mu} \d \mu_{\infty}$. Therefore, $\omega$ concentrates in $M_K$.
    \end{proof}
\end{proposition}

\begin{proposition}
    Let $\mu \in \mathscr{P}(X)$, assume that $\grad p \in C(\overline{\Omega})$ for $p = p[\mu]$ and define $\alpha:= \min_{x \in \overline{\Omega}} | \grad p(x)|$, $\beta := \max_{x \in \overline{\Omega}} | \grad p(x)|$. Then:
    \begin{enumerate}[label=(\roman*)]
        \item
        for $\gamma > 1$, $\mathrm{range}\left( \frac{\delta E[\mu]}{\delta \mu} \right) = \bigg[ - \frac{\gamma-1}{\gamma} \frac{\beta^{\frac{2\gamma}{\gamma-1}}}{\nu^{\frac{1}{\gamma-1}}}, \infty \bigg)$,
        \item
        for $0 < \gamma < 1$, $\mathrm{range}\left( \frac{\delta E[\mu]}{\delta \mu} \right) = \bigg( - \infty, \frac{1-\gamma}{\gamma} \frac{1}{\nu^{\frac{1}{\gamma-1}} \alpha^{\frac{2\gamma}{1-\gamma}}} \bigg]$,
        \item
        for $\gamma = 1$,
                \begin{itemize}
                \item
                if $\beta^2>\nu$, then $\mathrm{range}\left( \frac{\delta E[\mu]}{\delta \mu} \right) = (-\infty, +\infty)$,
                \item 
                if $\beta^2\leq\nu$, then $\mathrm{range}\left( \frac{\delta E[\mu]}{\delta \mu} \right) = [0, +\infty)$.
        \end{itemize}
    \end{enumerate}
    
    \begin{proof}
    Fix $y \in \R^d$ and consider
    \begin{align*}
        f_y (\C) := \frac{\nu}{\gamma} | \C |^{\gamma} - y \cdot \C y.
    \end{align*}
    If $\gamma > 1$, we have $\lim_{|\C| \rightarrow \infty} f_{y}(\C) = \infty$ and $f_y$ is convex. Therefore, $f_y$ attains its minimum. A simple computation shows that
    \begin{align*}
        \C_{\mathrm{min}} = \left( \frac{1}{\nu |y|^{2(\gamma-2)}} \right)^{\frac{1}{\gamma-1}} y \otimes y, \; \; f_y (\C_{\mathrm{min}}) = - \frac{\gamma-1}{\gamma} \frac{|y|^{\frac{2\gamma}{\gamma-1}}}{\nu^{\frac{1}{\gamma-1}}} < 0.
    \end{align*}
    Thus, $\mathrm{range}(f_y) = \bigg[ - \frac{\gamma-1}{\gamma} \frac{|y|^{\frac{2\gamma}{\gamma-1}}}{\nu^{\frac{1}{\gamma-1}}}, \infty \bigg)$ and since $|y|:=|\grad p(x)| \le \beta$ we conclude (i). To prove (ii), an analogous reasoning can be applied.

    For (iii), fix any $y^\perp\neq 0$ orthogonal to $y\in\R^d$ and set $\C^\perp_\lambda:=\lambda y^\perp \otimes y^\perp$ with $\lambda>0$. Then
    \[
       f_y(\C^\perp_\lambda) = \lambda\nu |y^\perp|^2 \to +\infty \quad\mbox{as } \lambda\to+\infty.
    \]
    Moreover, from the Cauchy-Schwartz inequality we have
    \[
       f_y(\C) \geq \left(\nu - |y|^2\right) |\C|
       \geq \left( \nu - \beta^2\right) |\C|,
    \]
    with equality for $\C_\lambda:=\lambda y\otimes y$ with $|y|=\beta$ and $\lambda>0$.
    Obviously, if $\nu - \beta^2 <0$,
    then
    $f_y(\C_\lambda) = \left( \nu - \beta^2\right) \lambda\beta^2 \to -\infty$ as $\lambda\to+\infty$.
    On the other hand, if $\nu - \beta^2 \geq 0$,
    then $f_y(\C_\lambda) \geq 0$ with equality for $\lambda:=0$.
    \end{proof}
\end{proposition}

In the subsequent discussion, the one-dimensional version of \eqref{eq:algebraic_equation_fr}, i.e., $\frac{\nu}{\gamma} C^{\gamma} - y^2 C - K = 0$ for $C\ge0$ and $y \in \R$, will play a pivotal role. Define $w := y^2 := p_x^2[\mu]$, we look for $C(w) \ge 0$ such that
\begin{align}
    \label{eq:algebraic_equation_fr_1D}
    \frac{\nu}{\gamma} C^{\gamma} - w C - K = 0, \quad \mathrm{for} \; w \in \R^+_0, C \in \R^+_0.
\end{align}
If such a function $C=C(w)$ exists, then the one-dimensional steady state solution of \eqref{eq:fr_gradientflow} reads
\begin{align}
    \label{eq:fr_stat_solution_1D}
    \mu_{\infty}(x,C) = \rho(x) \delta (C - C(p_x^2)),
\end{align}
with $\rho \in L^1(\Omega)$ and $p$ solves the nonlinear Poisson equation
\begin{align*}
    - \left( \left( r + \rho C(p_x^2) \right) p_x \right)_x = S,
\end{align*}
subject to zero flux boundary conditions.
Hence, our focus lies in determining $C(w)$, seeking explicit solutions where feasible or studying its properties when explicit solutions are not available.

\begin{proposition}\label{prop:11}
    There are either $0, 1$ or $2$ non-negative solutions $C = C(w)$ of \eqref{eq:algebraic_equation_fr_1D} for $w \ge 0$. If
    \begin{enumerate}
        \item \label{case1:fr_gamma>1} $\gamma > 1$ and
        \begin{enumerate}[label=(\roman*)]
            \item \label{case1.1:fr_Kpositive} $K > 0$ then there exists a unique positive solution $C = C(w)$ defined for all $w \ge 0$ with $\frac{\d C}{\d w}(w) > 0$ for all $w \ge 0$, $\lim_{w \rightarrow \infty} \frac{C(w)}{w^{\frac{1}{\gamma-1}}} = \left( \frac{\gamma}{\nu} \right)^{\frac{1}{\gamma-1}}$.
            \item \label{case1.2:fr_Kzero} $K = 0$ then $C_1(w) = 0$ and $C_2(w) = \left( \frac{\gamma w}{\nu} \right)^{\frac{1}{\gamma-1}}$ are the two solutions branches for $w \ge 0$.
            \item \label{case1.3:fr_Knegative} $K < 0$ then there are two solution branches $C_1(w)$ and $C_2(w)$ defined for $w \ge w_{\mathrm{min}, K} := \left( \frac{\gamma |K| \nu^{\frac{1}{\gamma-1}}}{|1 - \gamma|} \right)^{\frac{\gamma-1}{\gamma}}$ with $C_1(w) > C_2(w)$ for $w > w_{\mathrm{min}, K}$ and $C_1(w_{\mathrm{min}, K}) = C_2(w_{\mathrm{min}, K})$. Moreover, $\frac{\d C_1(w)}{\d w}(w) > 0$ and $\frac{\d C_2(w)}{\d w}(w) < 0$ for $w > w_{\mathrm{min}, K}$.
        \end{enumerate}
        \item \label{case2:fr_gamma=1} $\gamma = 1$ and
        \begin{enumerate}[label=(\roman*)]
            \item \label{case2.1:fr_Kpositive} $ K > 0$ then the only non-negative solution is $C(w) = \frac{K}{\nu - w}$ for $ 0 \le w < \nu$.
            \item \label{case2.2:fr_Kzero} $K = 0$ then $C (w) = 0$ for all $w \ge 0$ and $w = \nu$ solves \eqref{eq:algebraic_equation_fr_1D} for all $C \ge 0$.
            \item \label{case2.3:fr_Knegative} $K < 0$ then the only non-negative solution is $C(w) = \frac{|K|}{w - \nu}$ for $ w > \nu$.
        \end{enumerate}
        \item \label{case3:fr_gamma<1} $0 < \gamma < 1$ and
        \begin{enumerate}[label=(\roman*)]
            \item \label{case3.1:fr_Kpositive} $K > 0$ then there are two solution branches $C_1(w)$ and $C_2(w)$ defined for $0 < w \le w_{\mathrm{max}, K} := \left( \frac{1 - \gamma}{\gamma K \nu^{\frac{1}{\gamma-1}}} \right)^{\frac{1- \gamma}{\gamma}}$ with $C_1(w) < C_2(w)$ for $w < w_{\mathrm{max}, K}$ and $C_1(w_{\mathrm{max}, K}) = C_2(w_{\mathrm{max}, K})$. Moreover, $\frac{\d C_1(w)}{\d w}(w) > 0$ and $\frac{\d C_2(w)}{\d w}(w) < 0$ for $0 < w < w_{\mathrm{max}, K}$ and $C_1(0) = \left( \frac{\gamma K}{\nu} \right)^{\frac{1}{\gamma}}$ and $\lim_{w \rightarrow 0^+} C(w) = \infty$.
            \item \label{case3.2:fr_Kzero} $K = 0$ then $C_1(w) = 0$ and $C_2(w) = \left( \frac{\gamma w}{\nu} \right)^{\frac{1}{\gamma-1}}$ are the two solutions branches for $w > 0$ (clearly $C_1(0)=C_2(0) = 0$).
            \item \label{case3.3:fr_Knegative} $K < 0$ then there exists a unique positive solution $C = C(w)$ defined for all $w > 0$ with $\frac{\d C}{\d w}(w) < 0$ and $\lim_{w \rightarrow 0^+} C(w) = \infty $.
        \end{enumerate}
    \end{enumerate}
    \begin{proof}
        We prove case \ref{case1:fr_gamma>1}. Let us start with (i) taking $K > 0$. Notice that $C(0) = \left( \frac{\gamma K}{\nu} \right)^{\frac{1}{\gamma}} > 0$ and differentiating \eqref{eq:algebraic_equation_fr_1D}, we obtain
        \begin{align}
            \label{eq:algebraic_equation_deriv_fr_1D}
            C_w \left( \nu C^{\gamma-1} - w \right) = C.
        \end{align}
        It is evident that $C_w (0) > 0$, implying $C_w (w) > 0$ for $w \in [0, w_0)$, with $w_0 \in (0, \infty)$ sufficiently small. We claim that $w_0 = \infty$. Assuming for contradiction that $w_0 < \infty$, we have $C_w(w_0)=0$, leading to $C(w_0)=\infty$ from \eqref{eq:algebraic_equation_deriv_fr_1D}. Now, fix $\eps > 0$, and applying Young's inequality to the term $- w C$ yields
        \begin{align*}
            \frac{\nu}{\gamma} C^{\gamma} - \left( \frac{w}{\eps} \right)^{\gamma'} \frac{1}{\gamma'} - \frac{(\eps C)^{\gamma}}{\gamma} \le \frac{\nu}{\gamma} C^{\gamma} - w C = K, \quad \mathrm{with} \; \frac{1}{\gamma} + \frac{1}{\gamma'} = 1.
        \end{align*}
        Consequently,
        \begin{align*}
            \frac{1}{\gamma} \left( \nu - \eps^{\gamma} \right) C^{\gamma} \le K + \left( \frac{w}{\eps} \right)^{\gamma'} \frac{1}{\gamma'}
        \end{align*}
        and by choosing $\eps$ sufficiently small, we conclude that $C=C(w)$ is locally bounded on $[0, \infty]$. This contradicts $C(w_0) = \infty$ for $w_0<\infty$. Thus, $w_0 = \infty$ and $C=C(w)$ is strictly increasing on $[0, \infty)$. 
        This implies that $C$ is bounded below on any compact subset of the positive reals. Then, $\frac{\nu}{\gamma} C^{\gamma} = C w + K \ge \mathrm{const}(1 + w)$ and $C \ge \mathrm{const} (1 + w)^{\frac{1}{\gamma}}$, which also proves $\lim_{w \rightarrow \infty} C(w) = \infty$. Then, we can write $\frac{\nu}{\gamma} C^{\gamma-1} = w + \frac{K}{C}$ and taking the limit as $w \rightarrow \infty$ we conclude.

        The case (ii) is obvious.

        To prove (iii) we take $K < 0$ and define
        \begin{align*}
            f_w (C) = \frac{\nu}{\gamma} C^{\gamma} - w C.
        \end{align*}
        Since $\gamma>1$ we have $\lim_{C \rightarrow \infty} f_w (C) = \infty$ and $f_w$ is convex. Therefore, $f_w(C)$ attains its minimum. A simple computation shows
        \begin{align*}
            C_{\mathrm{min}} = \left( \frac{w}{\nu} \right)^{\frac{1}{\gamma-1}}, \; \; f_w (C_{\mathrm{min}}) = \left( \frac{1-\gamma}{\gamma} \right) \frac{w^{\frac{\gamma}{\gamma-1}}}{\nu^{\frac{1}{\gamma-1}}}.
        \end{align*}
        As $K < 0$, we  must require $| f_w(C_{\mathrm{min}}) | \ge | K |$, that is $w \ge w_{\mathrm{min}, K} := \left( \frac{\gamma |K| \nu^{\frac{1}{\gamma-1}}}{|1 - \gamma|} \right)^{\frac{\gamma-1}{\gamma}}$. Since $f_w(0)=0, f_w(C_{\mathrm{min}})<0$ and from convexity, we conclude that, for every $K$-level set we have two different points of intersection, so two branches $C_1(w)$ and $C_2(w)$. For one of the two branches, say $C_2(w)$, the term $\nu C_2^{\gamma-1}(w) - w < 0$ which implies $\frac{d C_2}{d w} < 0$ from \eqref{eq:algebraic_equation_deriv_fr_1D}. Same reasoning applies for $C_1(w)$.

        Proving \ref{case2:fr_gamma=1} is straightforward and \ref{case3:fr_gamma<1} follows the same steps of \ref{case1:fr_gamma>1}.
    \end{proof}
\end{proposition}

Going back to the multidimensional case, we write $\C = | \C | O^T \mathrm{diag}(\alpha) O$, where $O$ is an orthonormal matrix of eigenvectors of $\C$ and $\alpha \in \mathbb{S}^{d-1}_+ := \{ \omega \in \mathbb{S}^{d-1} : \omega = (\omega_1, \cdots, \omega_d) \; \mathrm{and} 
\; \omega_i \ge 0 \; \mathrm{for} \; i=1, \cdots, d \}$. Define $u:=|\C|$, $A:= O^T \mathrm{diag}(\alpha) O$ and $v:= y \cdot A y$ as new variables. Then, we can rewrite \eqref{eq:algebraic_equation_fr}
\begin{align}
    \label{eq:algebraic_equation_newvar_fr}
    \frac{\nu}{\gamma}u^{\gamma}-uv-K=0.
\end{align}
Note that with these new variables, we reduced \eqref{eq:algebraic_equation_fr} to the one-dimensional one \eqref{eq:algebraic_equation_fr_1D} studied above. In order to complete the construction of stationary solutions we (at first) take $\gamma>1$ and $K\ge0$. Later on we shall comment on the other cases. Again, the function $u=u(v)=u \left( \grad p \cdot A \grad p \right)$ is the solution $C = C(v)$ of \eqref{eq:algebraic_equation_newvar_fr}.

Further, let $X_1 := \overline{\Omega} \times \overline{\mathscr{S}^d_{1}}(\R)$, $\sigma_0 = \sigma_0(x, A) \in \mathscr{P}\left( X_1 \right)$, where $\overline{\mathscr{S}^d_{1}}(\R)$ is the $\left( \frac{d(d+1)}{2} - 1 \right)$-dimensional manifold of non-negative definite symmetric matrices with norm $1$, i.e., $\overline{\mathscr{S}^d_{1}}(\R) := \big{\{} A \in \overline{\mathscr{S}^d}(\R) : | A | = 1 \big{\}}$, and define 
\begin{align}
    \label{eq:sigma_stat}
    \sigma(x, A, r) := \sigma_0 (x, A) \delta \left( r - u \left( \grad p_{\infty} \cdot A \grad p_{\infty} \right) \right) \in \mathscr{P}\left( X_1 \times [0, \infty) \right).
\end{align}
Here, $u = u(v)$ is the solution of \eqref{eq:algebraic_equation_newvar_fr}.
\begin{proposition} \label{prop:fr_stationary}
    Let $Y := \overline{\mathscr{S}^d_{1}}(\R) \times [0, \infty)$ and define the (bijective) operator $T : Y \rightarrow \overline{\mathscr{S}^d}(\R)$ as $T(A, r) = r A$. Then
    \begin{align}
        \label{eq:fr_stationary_sol}
        \mu_{\infty}(x, \C) = (\mathrm{id}, T)_{\#} \sigma \in \mathscr{P}(X)
    \end{align}
    is a stationary solution of \eqref{eq:fr_gradientflow} if and only if $\mathbb{P}[\mu_{\infty}] \in L^1 (\Omega)$ and $p_{\infty} = p[\mu_{\infty}]$.
    \begin{proof}
        Let $\varphi(x, \C)$ be a test function and compute
        \begin{align*}
            \int_X \varphi \d \mu_{\infty} &= \int_X \varphi \d \left( (\mathrm{id}, T)_{\#} \sigma \right) \\
            &= \int_{\overline{\Omega}\times Y} \varphi(x, T(A, r)) \delta \left(r - u \left( \grad p_{\infty} \cdot A \grad p_{\infty} \right) \right) \d \sigma_0 \d r \\
            &= \int_{X_1} \varphi \left( x, u \left( \grad p_{\infty} \cdot A \grad p_{\infty} \right) A \right)\d \sigma_0.
        \end{align*}
        If $\varphi(x, \C) = \alpha(x) \beta(\C) \left( \frac{\nu}{\gamma} | \C |^{\gamma} - \grad p_{\infty} \cdot \C \grad p_{\infty} - K \right)$, where $\alpha, \beta$ are smooth test functions, we conclude that the above integral is zero by construction since $u$ is the solution of equation \eqref{eq:algebraic_equation_newvar_fr}.
    \end{proof}
\end{proposition}
We remark that $Y$ is a $\frac{d(d+1)}{2}$ dimensional manifold, isomorphic to $\overline{\mathscr{S}^d}(\R).$

At this point, we compute the permeability tensor
\begin{align*}
    \mathbb{P}[\mu_{\infty}](x) = \int_{\overline{\mathscr{S}^d_{1}}(\R)} u \left( \grad p_{\infty} \cdot A \grad p_{\infty} \right) A \sigma_0 \left( x, \d A \right).
\end{align*}

Let $U'(v) := u(v)$ with $U(0) = 0$ and define the functional $F : \overline{H}^1(\Omega) \rightarrow \R \cup \{ \infty \}$ by
\begin{eqnarray}
    \label{eq:functional_F}
        {\mathcal{F}[p] :=} \left\{ 
        \begin{array}{ll} \displaystyle
        \frac{r}{2} \int_{\Omega} | \grad p |^2 \d x + \frac{1}{2} \int_{X_1} U(\grad p \cdot A \grad p) \d \sigma_0 &\quad \mathrm{if} \int_{X_1} U( \grad p \cdot A \grad p) \d \sigma_0 < \infty \\ \displaystyle
        + \infty &\quad \mathrm{otherwise}.
        \end{array}
        \right.
\end{eqnarray}
We prove:
\begin{theorem}
    Let $\gamma > 1, K \ge 0, S \in L^2(\Omega)$ and $\int_{\overline{\mathscr{S}^d_1}(\R)} \sigma_0 ( \cdot, \d A) \in L^{\infty}(\Omega)$. Then, there exists a unique stationary state $\mu_{\infty}$ of \eqref{eq:fr_gradientflow} concentrating in $S_K$. $\mu_{\infty}$ is given by \eqref{eq:sigma_stat}-\eqref{eq:fr_stationary_sol}, where $p_{\infty}$ is the unique minimizer of the functional \eqref{eq:functional_F} in $\overline{H}^1(\Omega)$.
    \begin{proof}
        Note that (\ref{case1.1:fr_Kpositive}) of Proposition \ref{prop:11} implies the existence of positive constants $B_1, B_2$ such that
        \begin{align*}
            B_1 \left( 1+v^{\frac{1}{\gamma-1}} \right) \le u(v) \le B_2 \left( 1+v^{\frac{1}{\gamma-1}} \right), \; v \ge 0.
        \end{align*}
        Thus, integrating,
        \begin{align*}
            B_1 v \left( 1+ \frac{\gamma-1}{\gamma}v^{\frac{1}{\gamma-1}} \right) \le U(v) \le B_2 v \left( 1+\frac{\gamma-1}{\gamma}v^{\frac{1}{\gamma-1}} \right), \; v \ge 0.
        \end{align*}
        It is easy to conclude that, for some $B_3 > 0$
        \begin{align*}
            \int_{\Omega} \big{|} \mathbb{P}[\mu] \big{|} \d x \le B_3 \left( 1 + \int_{X_1} U(\grad p \cdot A \grad p) \d \sigma_0 \right)
        \end{align*}
        holds and $\mathbb{P}[\mu] \in L^1(\Omega)$ if $F[p] < \infty$. The same conclusion holds for case (\ref{case1.2:fr_Kzero}) of Proposition \ref{prop:11}. It is a simple exercise of variational calculus to show that critical points $p_{\infty}$ of $F$ are weak solutions of the Poisson equation
        \begin{align*}
            - \diver \left( \left( r \mathbb{I} + \mathbb{P}[\mu_{\infty}] \right) \grad p_{\infty} \right) &= S \quad \mathrm{in} \; \Omega \\
            n \cdot \left( r \mathbb{I} + \mathbb{P}[\mu_{\infty}] \right) \grad p_{\infty} &= 0 \quad \mathrm{on} \; \partial \Omega,
        \end{align*}
        (see also Theorem \ref{stationary_wasstype_gamma>1}).

        Define now $G(q) := \frac{1}{2} U \left( |q|^2 \right)$ for $q \in \R^d$ and compute the Hessian
        \begin{align*}
            D^2_q G(q) = 2 u'\left(|q|^2 \right) q \otimes q + u \left( |q|^2 \right) \mathbb{I}.
        \end{align*}
        Proposition \ref{prop:11}, (\ref{case1.1:fr_Kpositive}) and (\ref{case1.2:fr_Kzero}) give $D^2_q G(q) > 0$ on $\R^d$. The result then follows from classical variation theory \cite{evans}, taking into account that $D(F)$ is convex because of the convexity of $G(q)$ and that $C^1(\overline{\Omega}) \cap \big{\{} p : \int_{\Omega} p \d x = 0 \big{\}} \subseteq D(F)$.
    \end{proof}
\end{theorem}
We remark that the Theorem can easily be extended to the case of vanishing background permeability $r = 0$, if there exists $C_1 > 0$ such that
\begin{align}
    \label{eq:growth_sigma0}
    \int_{\overline{\mathscr{S}^d_1}(\R)} | z \cdot A z |^{\frac{\gamma}{1-\gamma}} \sigma_0 (x, \d A) \ge C_1 | z |^{\frac{2\gamma}{1-\gamma}} \quad \forall z \in \R^d, \; \mathrm{for} \; \mathrm{a.e.} \; x \in \Omega.
\end{align}
The uniqueness of $p_{\infty}$ in $W^{1, \frac{2\gamma}{1-\gamma}}(\Omega)$ follows.
Note that \eqref{eq:growth_sigma0} is satisfied, if, for example,
\begin{align}
    \label{eq:special_sigma0}
    d \sigma_0 = \rho(x) d x \times \left( V_{\#} \left( d \lambda \times d H \right) \right),
\end{align}
where $\rho \in L^1(\Omega)$ is uniformly positive, $V : \mathbb{S}^{d-1}_+ \times O(d) \rightarrow \overline{\mathscr{S}^d_1}(\R)$ with $V( \lambda, O) = O^T \Lambda O $ and $ \Lambda = \mathrm{diag}(\lambda)$, $d \lambda$ is a probability measure on $\mathbb{S}^{d-1}_+$, which charges a compact subset of the interior of $\mathbb{S}^{d-1}_+$ and $d H$ is a probability measure on $O(d)$ (real orthonormal matrices on $\R^d$).
Note that if \eqref{eq:special_sigma0} holds then the case (\ref{case1.2:fr_Kzero}) of Proposition \ref{prop:11} ($\gamma > 1, K = 0$) reduces exactly to the steady state case of Theorem \ref{stationary_wasstype_gamma>1} for the reduced Wasserstein flow. We also remark that case (\ref{case2.2:fr_Kzero}) of Proposition \ref{prop:11} ($\gamma=1, K = 0$) leads to the gradient constrained problem \eqref{eq:twisting}.
In all other cases of Proposition \ref{prop:11} at least one the following difficulties arises: either the branch $u = u(v)$ is represented by a strictly decreasing function of $v$ (which leads to an indefinite quadratic form in the Poisson energy $F[p]$) or the branch is not well-defined for all $v \ge 0$. This leads to constraints of either the form $ 0 \le \grad p(x) \cdot A \grad p(x) \le \alpha$ or $\grad p(x) \cdot A \grad p(x) \ge \alpha$ for some constant $\alpha > 0$.
Both of these issues will be addressed in future works.

\end{document}